\documentclass[10pt,sort&compress,3p]{elsarticle}

\usepackage{amsmath,amssymb,amsthm,mathrsfs,amsopn}
\usepackage{graphicx}
\usepackage{epstopdf}
\usepackage{subcaption}
\usepackage{booktabs}
\usepackage{array}
\usepackage{bm}
\usepackage{caption}
\usepackage{eucal}
\usepackage{tcolorbox} 
\tcbset{colback=blue!5!white,colframe=blue!50!black, left=0mm,right=0mm,top=0mm,bottom=0mm, box align=center, halign=center} 
\usepackage{tikz} 

\usetikzlibrary{automata}
\usepackage{enumitem}
\usepackage{fancyhdr} 
\usepackage{pdfpages} 
\usepackage{float}
\usepackage{multirow}
\usepackage{tabularx}
\usepackage{longtable}
\usepackage{booktabs}
\usepackage{diagbox}
\usepackage{rotating}
\usepackage{framed} 
\usepackage{multicol} 
\usepackage{nomencl}
\begin{document}

\begin{frontmatter}

\title{Integrated intelligent Jaya Runge-Kutta method for solving Falkner-Skan equations for Various Wedge Angles}
\author[add3]{Hongwei Guo\fnref{email1}}
\fntext[email1]{email: ghway0723{@}gmail.com}
\author[add3,add4,add5]{Xiaoying Zhuang\fnref{email2}}
\fntext[email2]{email: zhuang{@}ikm.uni-hannover.de}
\author[add3]{Xiaoyu Meng \fnref{email3}}
\fntext[email3]{email: mxy196711@gmail.com}
\author[add1,add2]{Timon Rabczuk\corref{mycorrespondingauthor}}
\cortext[mycorrespondingauthor]{Corresponding author}
\ead{timon.rabczuk@tdtu.edu.vn}

\address[add1]{Division of Computational Mechanics,\\ Ton Duc Thang University,\\ Ho Chi Minh City, Vietnam}
\address[add2]{Faculty of Civil Engineering,\\ Ton Duc Thang University,\\ Ho Chi Minh City, Vietnam}
\address[add3]{Chair of Computational Science and Simulation Technology,\\ Faculty of Mathematics and Physics,\\ Leibniz Universität Hannover,\\ Hannover, Germany}
\address[add4]{Department of Geotechnical Engineering,\\ Tongji University,\\ Shanghai, China.}
\address[add5]{Key Laboratory of Geotechnical and Underground Engineering of Ministry of Education,\\ Tongji University,\\ Shanghai, China.}
\begin{abstract}
In this work, the hybrid intelligent computing method, which combines efficient Jaya algorithm with classical Runge-Kutta method is applied to solve the Falkner-Skan equations with various wedge angles, which is the fundamental equation for a variety of computational fluid mechanical problems. With some coordinate transformation, the Falkner-Skan boundary layer problem is then converted into a free boundary problem defined on a finite interval. Then using higher order reduction strategies, the whole problem can be boiled down to a solving of coupled differential equations with prescribed initial and boundary conditions. The hybrid Jaya Runge-Kutta method is found to yield stable and accurate results and able to extract those unknown parameters. The sensitivity of classical shooting method to the guess of initial values can be easily overcome by an integrated robust optimization method. In addition, the Jaya algorithm, without the need for tuning the algorithm-specific parameters, is proved to be effective and stable for minimizing the fitness function in application. By comparing the solutions using the Jaya method with PSO (particle swarm optimization), Genetic algorithm (GA), Hyperband, and the classical analytical methods, the hybrid Jaya Runge-Kutta method yields more stable and accurate results, which shows great potential for solving more complicated multi-field and multiphase flow problems.

\end{abstract}
\end{frontmatter}

\section{Introduction}
The boundary layer flow problems over a stretching or shrinking surface have significant applications in many industrial and technological fields, such as cooling of an infinite metallic plate in a cooling bath, an aerodynamic extrusion of plastic sheets, and extraction of polymer and rubber and glass-fiber. The Falkner-Skan problem is a notable similarity solution to the steady two-dimensional laminar boundary layer equations, including the Blasius and stagnation point solutions.  

Historically, the Falkner-Skan equation was initially developed by Falkner and Skan in 1931 \cite{falkneb1931lxxxv}, and it plays an essential role in the fluid mechanics of boundary viscous flow. It is derived from the two-dimensional incompressible Navier-Stokes equations. It exists in many forms with varying values of $\beta_{0}$ and $\beta$, and its solution describes the forms on a wedge for the steady two-dimensional laminar boundary layer. In the past years, there are many investigations for the Falkner-Skan equation with numerical and analytical methods. The existence and uniqueness results were found by Rosenhead (1963), Weyl (1942), Hartman (1972), and Tam (1970) \cite{duque2011numerical}. The shooting and invariant imbedding is the earliest computational method, which is proposed by Hartree, Smith, Cebeci and Keller, Na \cite{asaithambi1998finite}. There are other numerical methods such as finite-element method by Asaithambi \cite{asaithambi2005solution}, the modification of the classical Newtons method by Zhang and Chen, etc.\cite{zhang2013particle}.

Because of limitations in the accuracy and efficiency of those classic and numerical methods in solving Falkner-Skan equations, the heuristic algorithms have been successfully applied to solve the nonlinear problems in engineering and applied science \cite{malik2015numerical,ullah2018evolutionary}. The heuristic algorithms are used to find a good solution effective and quickly, which are essentially the trial and error methods. Many researchers have devoted their attentions to study the heuristic algorithms and have thus developed a series of heuristic algorithms. The genetic algorithm(GA) is proposed by Vakil-Baghmisheh et al \cite{vakil2008crack} and Liu and Jiao\cite{liu2011application}. The particle swarm optimization(PSO) is found by Yu and Chen\cite{yu2010bridge}, Nanda,Maity,and Maiti\cite{nanda2014crack}. Zare Hosseinzadeh et al. proposed cuckoo search (CS) algorithm \cite{hosseinzadeh2014flexibility}. The evolutionary intelligence algorithms are proposed by Bagheri, Razeghi, and Ghodrati Amiri \cite{bagheri2012detection}. Recently, Li et al. \cite{li2017hyperband} developed a hyperband optimization method based on a novel bandit-based approach. The evolutionary computing algorithms have been proved to be effective to solve PDEs. Tsoulos and Lagaris \cite{tsoulos2006solving}, Saini et al. \cite{saini2013genetic} et al. applied genetic algorithm for solving partial differential equations. Babaei \cite{babaei2013general}, Yadav et al. \cite{yadav2017efficient} applied particle swarm method for solving PDEs. Also, the deep neural networks based energy method and collocation method have been applied to solving partial differential equations \cite{guo2019deep,lin2020deep,nguyen2019deep,anitescu2019artificial,samaniego2020energy}. However, the major drawback of those previous soft computing techniques, as well as the meta-heuristic, and evolutionary optimizations is the possibility of converging to solutions that may not be optimal but instead are trapped at a local optimal value \cite{wadood2019application}.
   
   Until recently, an advanced optimization algorithm called Jaya (a Sanskrit word meaning victory) was proposed by Venkata Rao in 2016 \cite{rao2016Jaya}, which is easy in implementation and does not need any algorithm-specific control parameters. The Jaya algorithm can solve both constrained and unconstrained optimization problems, and it is proved to converge to the global optimum values requiring less numbers of iterations and computational time  \cite{wadood2019application}. From the results we can prove that the Jaya algorithm is more effective, accurate and stable than other optimization algorithms.

 In this work, we first describe the mathematic model for Falkner-Skan problem and apply a transformation and order reduction to Falkner-Skan equation to obtain a system of coupled differential equation with prescribed boundary conditions. Then the optimization algorithms involved in the application are brief presented. The Jaya method is integrated with the classical Runge-Kutta method and used to solve the Falkner-Skan equation with various wedge angles. Also, as a comparison, the newly proposed hyperband method is also firstly applied to solve partial differential equations. The Runge-Kutta method combined with the Jaya algorithm is used to solve the Falkner-Skan equation in following manner:

 1. Using a coordinate transformation the semi-inifinite domain of Falkner-Skan equation can be changed into the finite unit interval.

 2. Using the Runge-Kutta method the solving of the Falkner-Skan equation can be converted into finding the initial values $\alpha$ and ${ \eta  }_{ \infty  }$.

 3. Based on the fitness function $min\quad F\left( \alpha ,{ \eta  }_{ \infty  } \right) = \parallel { RK }_{ 2 }\left( \alpha ,{ \eta  }_{ \infty  } \right) -1;{ RK }_{ 3 }\left( \alpha ,{ \eta  }_{ \infty  } \right)   \parallel $, the optimization algorithm find the optimal values of $\alpha$ and ${ \eta  }_{ \infty  }$.

 4.  The optimal initial values of $\alpha$ and ${ \eta  }_{ \infty  }$ are used for finding the solution of the Falkner-Skan equation with the multistep method.

\section{Mathematical models for Falkner-Skan equation}
In the case of the steady state boundary layer flow of an incompressible viscous fluid over a wedge, the equations of continuity in two-dimensional motions and their boundary conditions \cite{schlichting2016boundary} can be referred as:
The boundary layer equation for the steady flow of an incompressible viscous fluid over a wedge in two-dimensional is defined as:
\begin{equation}
\label{3.1}
\frac { \partial v }{ \partial x } +\frac { \partial w }{ \partial y } =0
\end{equation}
\begin{equation}
\label{3.2}
v\frac { \partial v }{ \partial x } +w\frac { \partial w }{ \partial y } =U\frac { dU }{ dx } +\mu \frac { { \partial  }^{ 2 }v }{ \partial { y }^{ 2 } } ,
\end{equation}
where $v$ and $w$ are velocity components in $x$ and $y$ direction of the fluid flow, $\mu$ is the viscosity, $U\left(x\right)$ the velocity at the edge of the boundary layer. Assuming the velocity of the ambient
flow and moving wedge flow is a power law free stream velocity, $U\left(x\right)=U_{\infty} \left(\frac{x}{L}\right)^m$, where $U_{\infty}$ is uniform free stream velocity, L is the length of wedge, $x$ is measured from the top of the wedge and m is the Falkner-Skan power-law parameter.

The relevant boundary conditions are given by:
\begin{equation}
\label{3.3}
\begin{aligned}
&v\left( x,0 \right) =w\left( x,0 \right) =0\\
&v\left( x,y \right)   \rightarrow   U\left( x \right) \;\;  as\;\; y \rightarrow  \infty 
\end{aligned}
\end{equation}

The continuity Equation \ref{3.1} is automatically satisfied by the stream function $\psi \left(x,y \right)$, thus obtaining:
\begin{equation}
\label{3.4}
 v=\frac { \partial \psi  }{ \partial y } \;\;    \textrm{and}\;\;   w=-\frac { \partial \psi  }{ \partial x } 
\end{equation}
and from momentum Equation \ref{3.2} yields
\begin{equation}
\label{3.5}
\frac { \partial \psi  }{ \partial y } \frac { { \partial  }^{ 2 }\psi  }{ \partial y\partial x } -\frac { \partial \psi  }{ \partial x } \frac { { \partial  }^{ 2 }\psi  }{ \partial { y }^{ 2 } } =U\frac { dU }{ dx } +\mu \frac { { \partial  }^{ 3 }\psi  }{ \partial { y }^{ 3 } } 
\end{equation}
Using the similarity transformation, two dimensionless variables can be yielded:
\begin{equation}
\label{3.6}
\begin{aligned}
&f\left( \eta  \right) =\sqrt { \frac { 1+m }{ 2 } \frac { { L }^{ m } }{ \mu { U }_{ \infty  } }  } \left[ \frac { \psi  }{ { x }^{ \frac { 1+m }{ 2 }  } }  \right] ,\\
&\eta =\sqrt { \frac { 1+m }{ 2 } \frac { { U }_{ \infty  } }{ \mu { L }^{ m } }  } \left[ \frac { y }{ { x }^{ \frac { 1-m }{ 2 }  } }  \right] ,
\end{aligned}
\end{equation}
, with $f$ a dimensionless stream function and $\eta$ a dimensionless distance. 
Then the partial differential equation can be reduced to a third order nonlinear ordinary differential equation:
\begin{equation}\label{fs}
\frac { { d }^{ 3 }f\left( \eta  \right)  }{ d{ \eta  }^{ 3 } } +{ \beta  }_{ 0 }f\left( \eta  \right) \frac { { d }^{ 2 }f\left( \eta  \right)  }{ d{ \eta  }^{ 2 } } +\beta \left[ 1-{ \left( \frac { df\left( \eta  \right)  }{ d\eta  }  \right)  }^{ 2 } \right] =0,    \eta  \in  \left[ 0,\infty  \right) 
\end{equation}
The associated boundary conditions are given by:
\begin{equation}\label{bd1}
\begin{aligned}
&f=0,\quad at  \quad  \eta =0,\\
&\frac { df }{ d\eta  } =0\quad at  \quad \eta =0. 
\end{aligned}
\end{equation}
and 
\begin{equation}
\label{asymp}
\frac { df }{ d\eta  } =1,\quad as  \quad\eta \rightarrow\infty.
\end{equation}
 where ${ \beta  }_{ 0 }$ and ${ \beta  }$ are constants, $\beta \pi$ is the angle of the wedge. $\beta$ is related to the Falkner-Skan power-law parameter $m$ by $\beta=\frac{2m}{1+m}$.

\begin{figure}[H]
	\captionsetup{width=0.9\columnwidth}
	\centering
	\begin{subfigure}[b]{6.0cm}
		\centering\includegraphics[height=5.5cm,width=7.0cm]{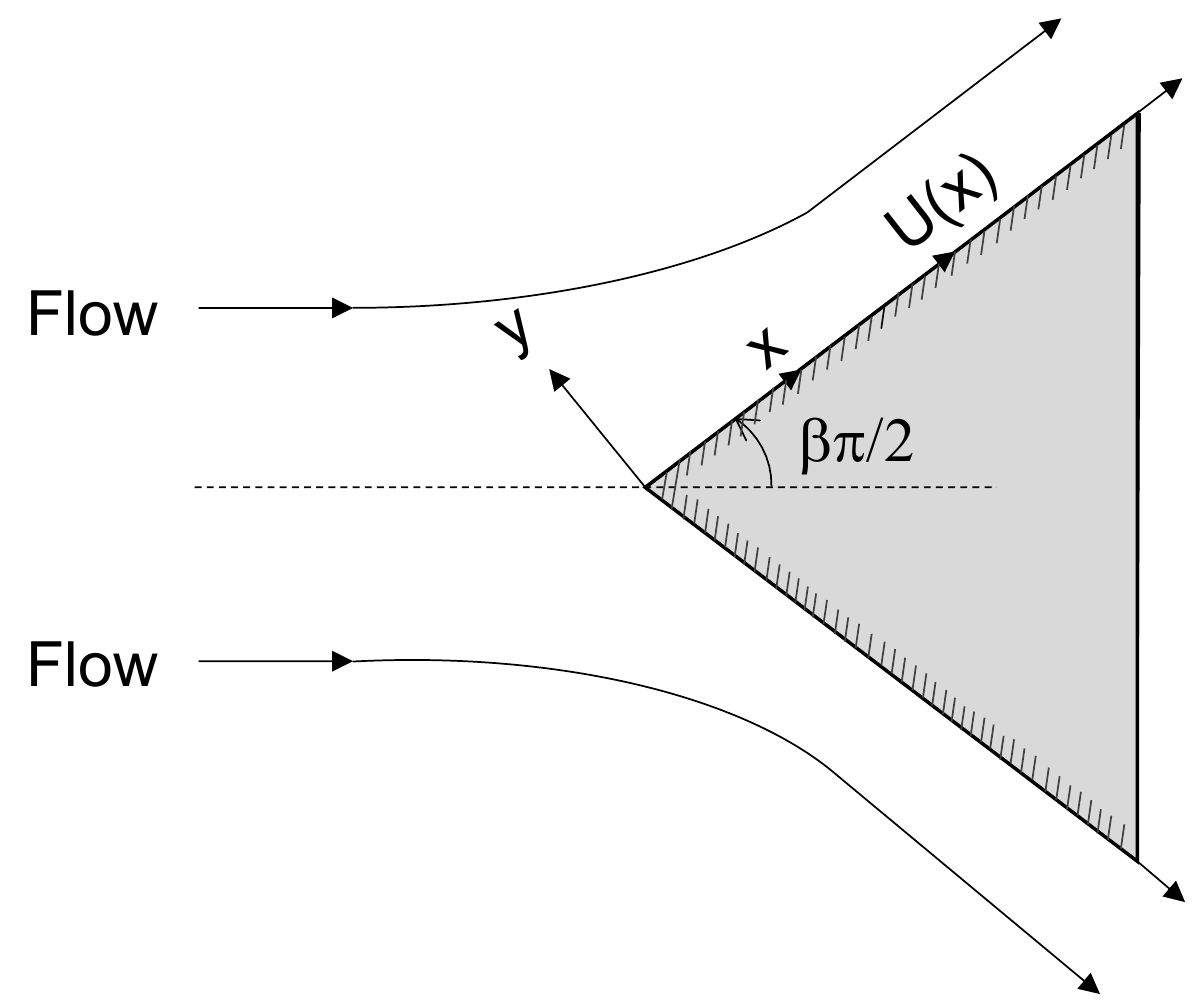}   
		\caption{}\label{flow1}
	\end{subfigure}%
	\hspace{1cm}
    \begin{subfigure}[b]{6.0cm}
	\centering\includegraphics[height=5.5cm,width=7.0cm]{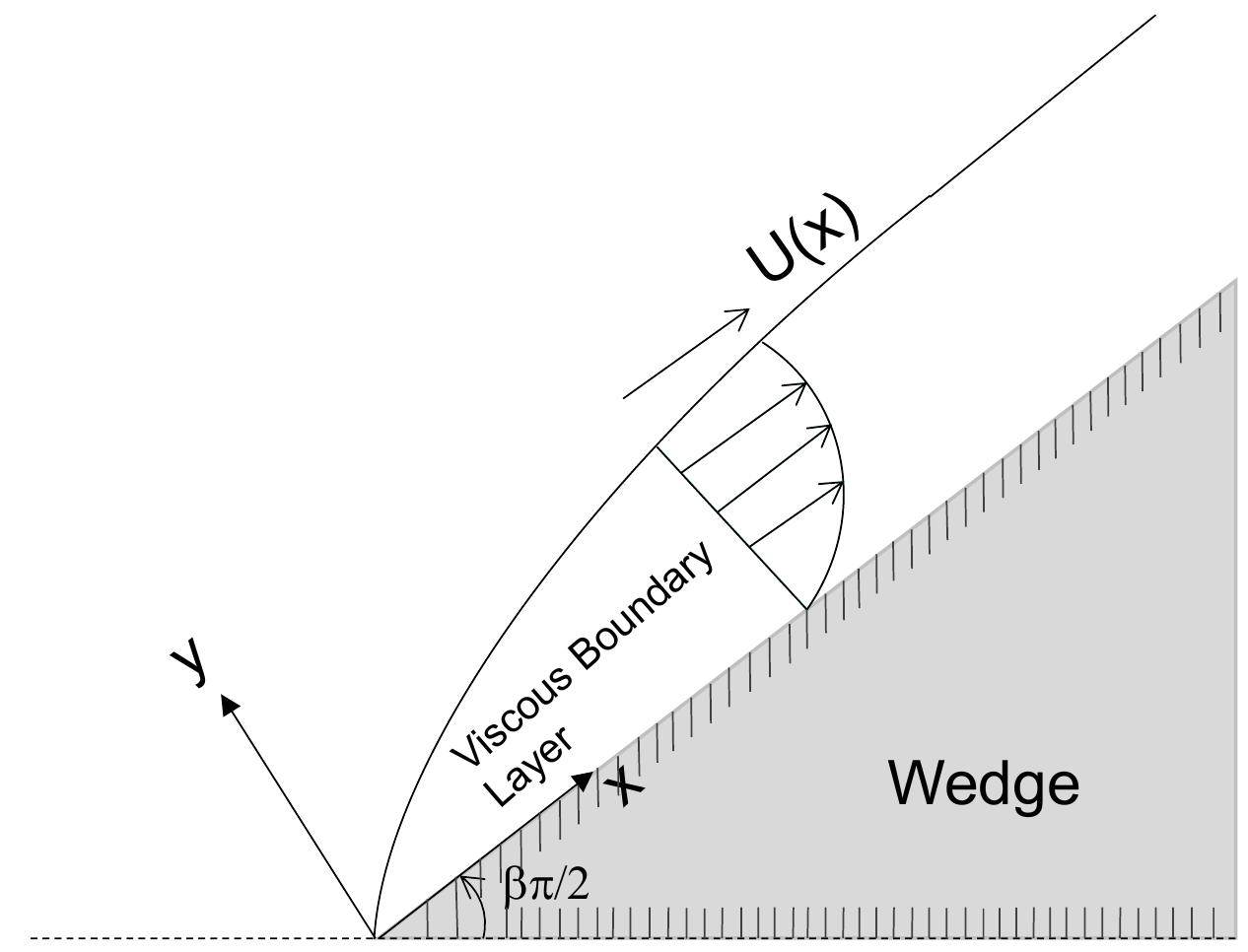}   
	\caption{}\label{flow2}
    \end{subfigure}%
    \newline
    \vspace{0.5cm}
    \begin{subfigure}[b]{6.0cm}
	\centering\includegraphics[height=5.5cm,width=7.0cm]{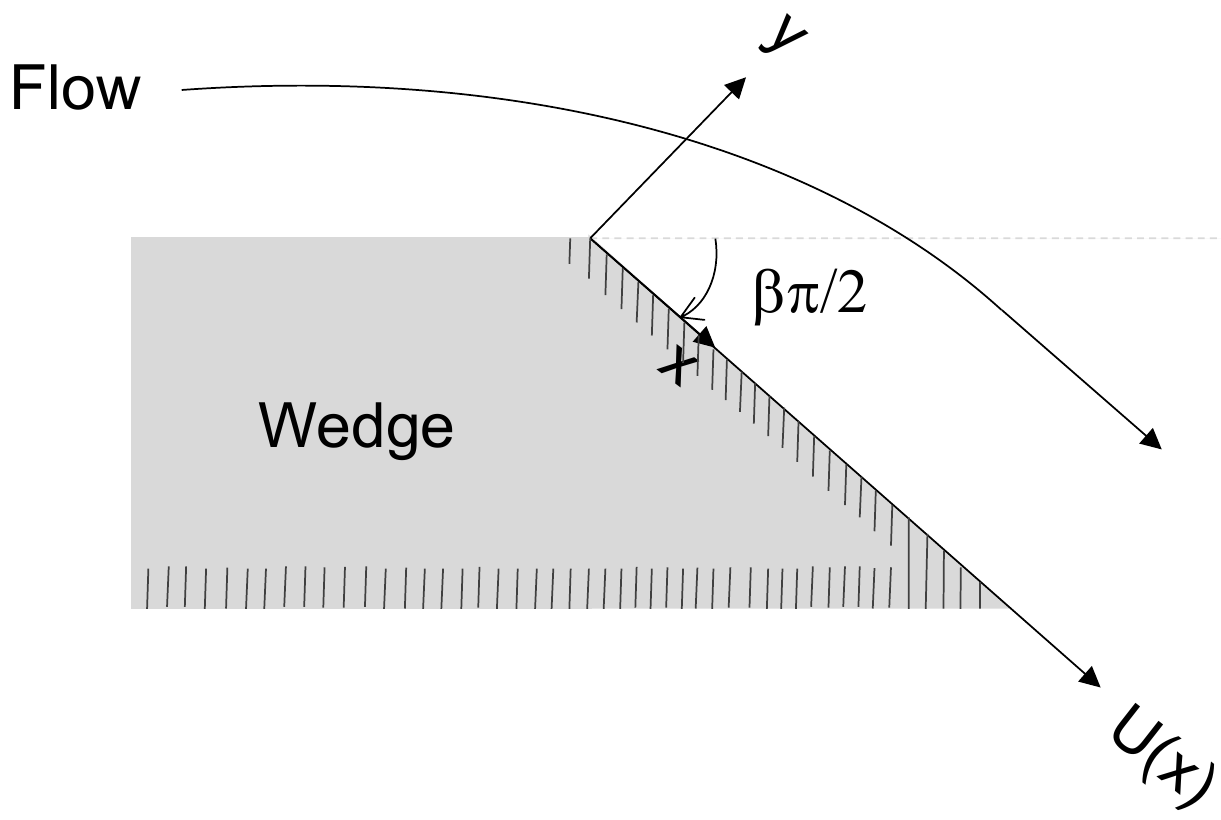}   
	\caption{}\label{flow3}
    \end{subfigure}%
	\hspace{1cm}
    \begin{subfigure}[b]{6.0cm}
	\centering\includegraphics[height=5.5cm,width=7.0cm]{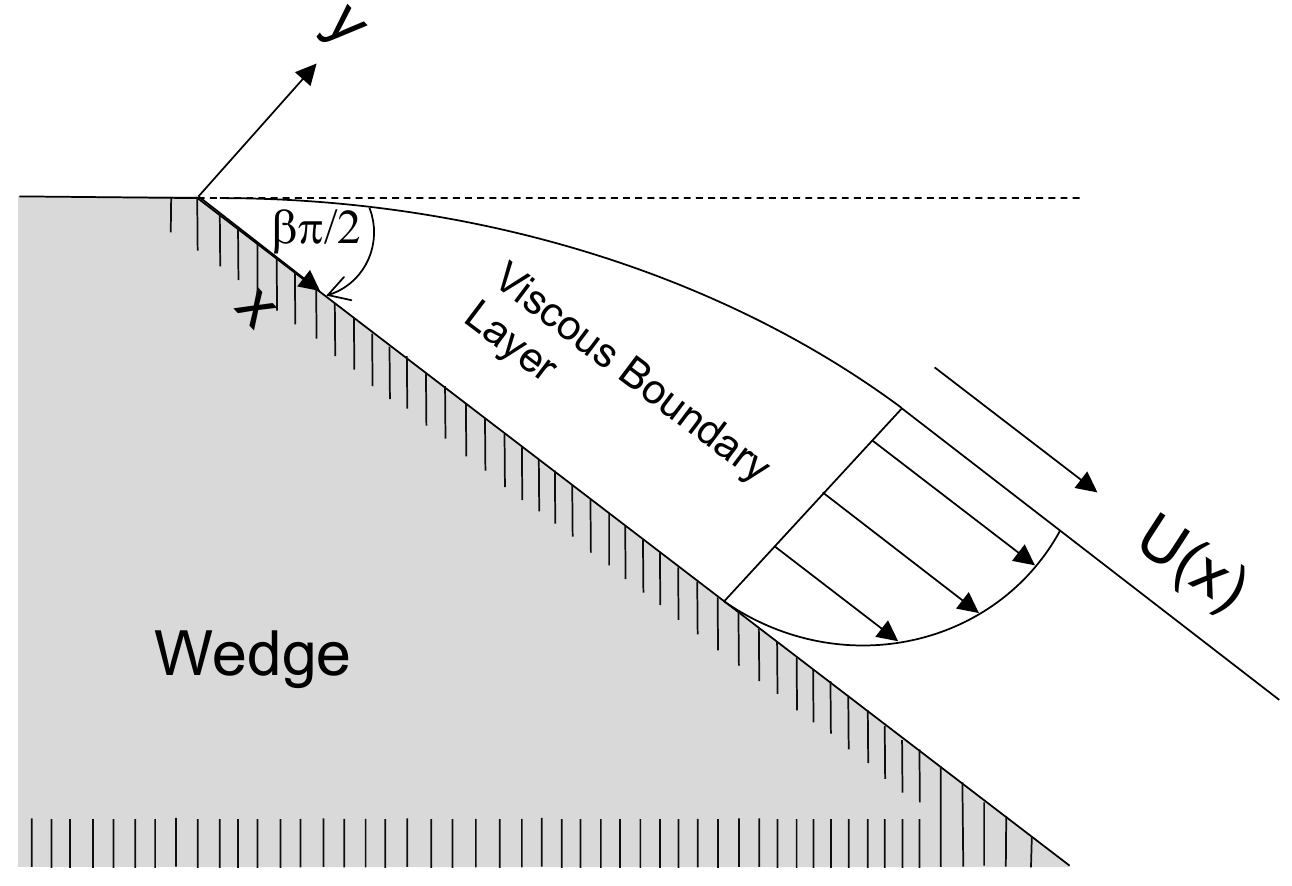}   
	\caption{}\label{flow4}
    \end{subfigure}%
	\caption{Geometry of Falkner-Skan flow over a wedge. a) Domain for the potential flow solution for $\beta>0$ with nonlinear velocity $U(x)$. b) Domain for the Falkner-Skan boundary layer. c) and d) are $\beta<0$ cases.}
	\label{flux_1}
\end{figure}
  
The Geometry of Falkner-Skan flow over a wedge is shown in Figure~\ref{flux_1}, which describes the potential flow and thin viscous boundary layer regions. The Falkner–Skan equation~\ref{fs} governing this phenomenon describes a nonlinear, third-order, free boundary value problems, where $f^{'}(\eta)$ defines the dimensionless velocity component indirection and $f^{''}(\eta)$ the dimensionless shear stress in the boundary layer. Moreover, the solution of this equation is sought in a semi-infinite domain and one of the boundary conditions is asymptotically assigned on the first derivative at infinity. Thus, it is not easily to construct a closed-form solution for this two-point boundary value problem. 
 
 \section{Transformation and order reduction}
 First of all, we replace the boundary condition at $\infty$ in Equation~\ref{asymp} with a free boundary condition as Asaithambi did in \cite{asaithambi2005solution}:
 \begin{equation}\label{bd2}
\frac { df }{ d\eta  } =1 , \;\; \eta =n_{\infty},
\end{equation}
where $n_{\infty}$ is an unknown truncated boundary (a free boundary). The whole problem is then converted into a free boundary problem defined on a finite interval, where the “sufficiently large” $\eta_{\infty }$ is determined as part of the solution.
 
To ensure the impose of asymptotic boundary condition, additionally boundary conditions need to be added for the semi-infinite physical domain \cite{asaithambi2005solution}:
\begin{equation}\label{bd3}
\frac { { d }^{ 2 }f }{ d{ \eta  }^{ 2 } } =0 , \;\; \eta =n_{\infty}.
\end{equation}

 If a shooting algorithm is imposed to solve this free boundary problems governed by Equation ~\ref{fs} subjected to boundary conditions Equations~\ref{bd1}, \ref{bd2}, and \ref{bd3}, an initial value condition is added at $\eta =0$:
\begin{equation}
\frac { { d }^{ 2 }f }{ d{ \eta  }^{ 2 } } =\alpha , \;\; \eta =0.
\end{equation}
It is obvious that different values of $\alpha$ will lead to different values of $\frac { df }{ d\eta  } $ as $\eta$ approaches infinity. The $\alpha$ value is generally used to characterize solutions of Falkner–Skan equation with different $\beta0$ and $\beta$, and it will be used for model evaluation. It is well known that the shooting method can fail to converge for problems whose solutions are very sensitive to initial conditions. Many strategies have been proposed to tackle this problem \cite{holsapple2004new}. A new, simple and straightforward approach has been proposed to the whole scheme, which deals with the unknown initial boundary conditions and the free boundary with a robust optimization method.

For simplicity, a coordinate transformation is first placed: 
\begin{equation}
\xi =\frac { \eta  }{ { \eta  }_{ \infty  } },
\end{equation}
which transforms the physical domain from $[0,n_{\infty})$ to $[0,1]$, and yields Falkner-Skan equation:
\begin{equation}
\frac { 1 }{ { \eta  }_{ \infty  }^{ 3 } } \frac { { d }^{ 3 }f }{ d{ \xi  }^{ 3 } } +\frac { 1 }{ { \eta  }_{ \infty  }^{ 2 } } { \beta  }_{ 0 }f\frac { { d }^{ 2 }f }{ d{ \xi  }^{ 2 } } +\beta \left[ 1-\frac { 1 }{ { \eta  }_{ \infty  }^{ 2 } } { \left( \frac { df }{ d\xi  }  \right)  }^{ 2 } \right] =0, 0<\xi <1 
\end{equation}
and the boundary conditions are changed correspondingly:
\begin{equation}
\begin{aligned}
&f=0,  \xi =0,\\
& \frac { df }{ d\xi  } =0,  \xi =0,\\
& \frac { df }{ d\xi  } ={ \eta  }_{ \infty  }, \xi =1,\\
& \frac { { d }^{ 2 }f }{ d{ \xi  }^{ 2 } } =0,  \xi =1,\\
& \frac { { d }^{ 2 }f }{ d{ \xi  }^{ 2 } } ={ \eta  }_{ \infty  }^{ 2 }\alpha ,  \xi =0.\\ 
\end{aligned}
\end{equation}

 In order to solve the Falkner-Skan equation and remove the $n_{\infty}^2$ and $n_{\infty}^3$ from above equations, the reduction of order technique is applied:
\begin{equation}
\begin{aligned}
&{ y }_{ 1 }=f,\\
& { y }_{ 2 }=\frac { 1 }{ { \eta  }_{ \infty  } } \frac { df }{ d\xi  } ,\\ 
&{ y }_{ 3 }=\frac { 1 }{ { \eta  }_{ \infty  }^{ 2 } } \frac { { d }^{ 2 }f }{ d{ \xi  }^{ 2 } } ,\\
& { y }_{ 4 }={ \eta  }_{ \infty  }
\end{aligned}
\end{equation}
 and the Falkner-Skan equation becomes
\begin{equation}\label{newfs}
{\begin{pmatrix}
	f_{1}\\ 
	f_{2}\\ 
	f_{3}\\ 
	f_{4}
	\end{pmatrix}}'=\begin{pmatrix}
{ f }_{ 2 }{ f }_{ 4 }\\ 
{ f }_{ 3 }{ f }_{ 4 }\\ 
-\beta _{0}{ f }_{ 1 }{ f }_{ 3 }{ f }_{ 4 }-\beta \left ( 1-{ f }_{ 2 }^{ 2 } \right ){ f }_{ 4 }\\ 
0
\end{pmatrix}
\end{equation}
 where $\eta \in [0,1] $ and the initial conditions conditions are:
\begin{equation}\label{newfsbd1}
\begin{aligned}
&{ f }_{ 1 }\left( 0 \right) =0,\\ 
&{ f }_{ 2 }\left( 0 \right) =0\\
& { f }_{ 3 }\left( 0 \right) =\alpha \\ 
&{ f }_{ 4 }\left( 0 \right) ={ \eta  }_{ \infty  }.
\end{aligned}
\end{equation}
the boundary conditions conditions are:
\begin{equation}\label{newfsbd2}
\begin{aligned}
&{ f }_{ 2 }\left( 1 \right) =1,\\
&{ f }_{ 3 }\left( 1 \right) =0,\\ 
&{ f }_{ 4 }\left( 1 \right) ={ \eta  }_{ \infty  }.
\end{aligned}
\end{equation}

\section{The Optimization Algorithm}
In this part, we will briefly introduce some classical heuristic algorithms and newly developed optimization algorithms, which will be used in the numerical example comparisons.
\subsection{PSO}
Inspired by social intelligent behaviors of the  gregarious birds the Particle swarm optimization (PSO) is proposed by Kennedy and Eberhar t\cite{kennedy1995particle}. PSO is a heuristic searching process. The search space of the optimization problem is similar to the flying space of birds, and the search process of the optimal solution can be described by the process of birds searching for food. In the process,  the local optimum position $P_{best}$ and the global optimal position $g_{best}$ in the population are obtained, with the help of current optimal particles the velocity of each particle will be adjusted to search for the optimal solution. SPO is a classical swarm intelligence approach. SPO has many advantages, its calculation method is relatively simple, and the method needs fewer control parameters, it can obtain a accurate optimal solution. But PSO also has a problem, due to the lack of individual diversity, it is easy to fall into local optimization with the iteration of the search process.\\
Assuming that a swarm has P particles, the position vector and the velocity vector at t iteration are ${ X }_{ i }^{ t }$ and ${ V }_{ i }^{ t }$, the vector is updated by following equation:
\begin{equation}
\begin{aligned}
\label{PSO}
&{ V }_{ ij }^{ t+1 }=w{ V }_{ ij }^{ t }+{ c }_{ 1 }{ r }_{ 1 }^{ t }\left( { pbest }_{ ij }-{ X }_{ ij }^{ t } \right) +{ c }_{ 2 }{ r }_{ 2 }^{ t }\left( { gbest }_{ j }-{ X }_{ ij }^{ t } \right) \\
&{ X }_{ ij }^{ t+1 }={ X }_{ ij }^{ t }+{ X }_{ ij }^{ t+1 }
\end{aligned}
\end{equation}  
where i=(1,2,...P) is the particle, and j=(1,2,...n) is the dimension, w is the inertia weight constant, which is a positive. The process of the PSO \cite{le2019hybrid} is presented in Fig.~\ref{PSO}.
 \begin{figure}[H]
	\captionsetup{width=0.9\columnwidth}
	\centering\includegraphics[height=14cm,width=14.0cm]{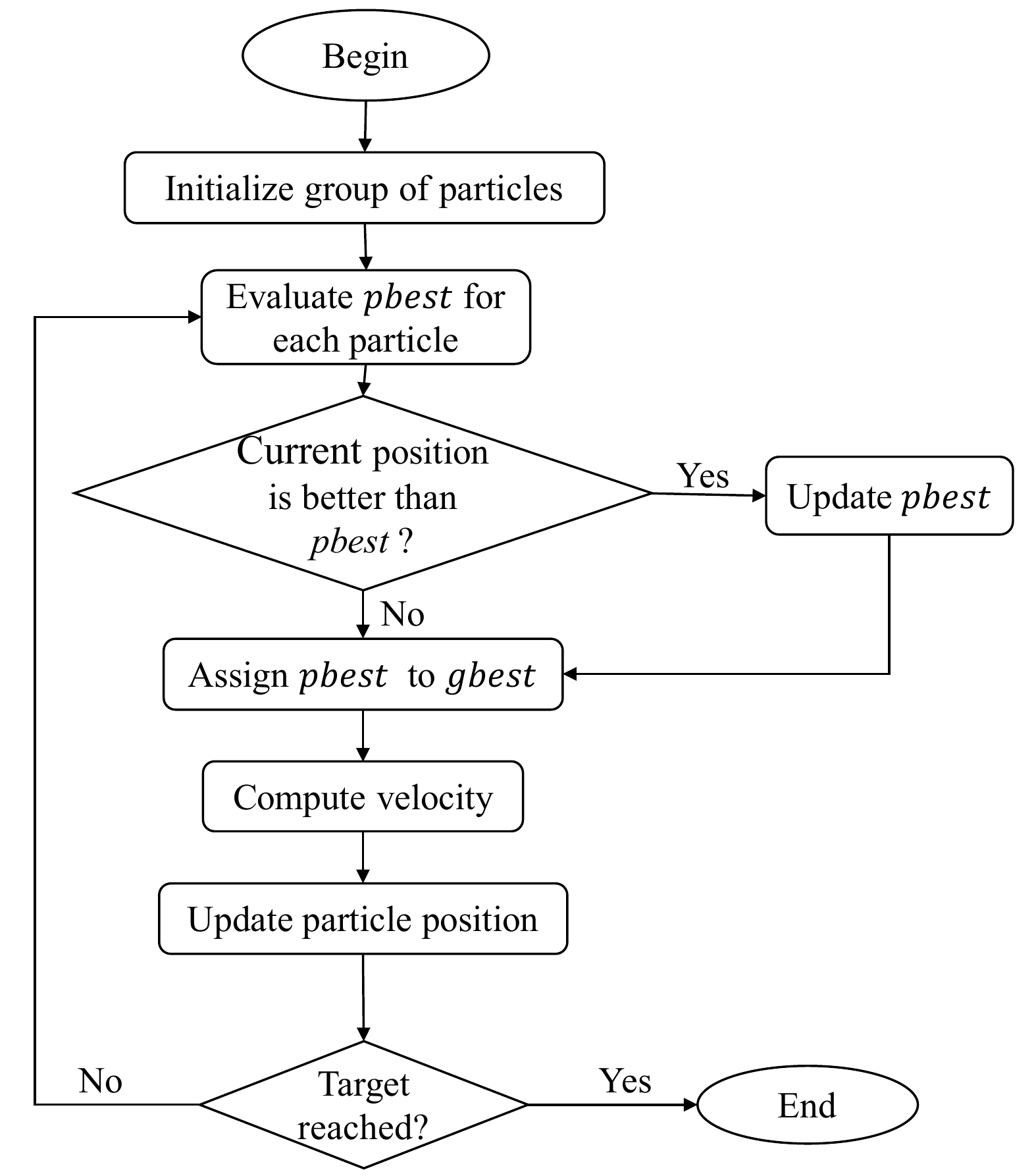} 
	\caption{ The Particle Swarm Optimization}
	\label{PSO}
\end{figure}

\subsection{Hyperband Algorithm}
Hyperband is a SuccessiveHalving Algorithm. The process can be described as follow. In the search space, we sample $n$ hyper-parameter sets randomly, after $m$ interactions the validation loss is evaluated, then we discard the lowest performers of the parameters and keep a half. in the next step, we use the good parameters to run for $m$ interactions, and discard a half. The process will be repeated until we have only one model. It's a variation of random search, and it can assign adaptive pre-defined parameters to a randomly sampled configuration automatically. In the process of optimization, Hyperband can find the best time allocation for each configuration. Compared with Bayesian Optimization methods on several hyperparameter optimization problems, the Hyperband can find the optimal  solution with faster speed, the performance of the optimization is good, and it can parallelize the operating parameters. The process of Hyperband Algorithm \cite{kainz2019efficient} is presented in Table.~\ref{Hyperband}.
\begin{table}[H]
	\captionsetup{width=0.85\columnwidth}
	\caption{\textbf {Hyperband ALgorithm}} 
	\vspace{-0.3cm}
	\centering
    \resizebox{0.85\columnwidth}{!}{%
        \begin{tabular}{l} 
            \toprule 
             \quad \textbf {input : R, $\eta$}\\ 
            1  \textbf {initialization}: ${ s }_{ max }=\left\lfloor { log }_{ \eta  }R \right\rfloor $ and $B=\left( { s }_{ max }+1 \right) R$;\\
            2 \quad \textbf {for} $s\in \left\{ { s }_{ max }\quad ,\quad { s }_{ max }-1,...,\quad 0 \right\}$ \quad \textbf {do}\\
            3 \quad \quad \quad  $n=\left\lceil \frac { B }{ R } \frac { { \eta  }^{ s } }{ s+1 }  \right\rceil $,  $r=R{ \eta  }^{ -s }$;\\
            4 \quad \quad \quad  ${ \Lambda  }_{ s }$= get$\_$hyperparameter$\_ $configuration(n);\\
            5 \quad \quad \quad  \textbf {for}  $i\in \left\{ 0,...,s \right\}$ \textbf {do}\\
            6 \quad \quad \quad \quad \quad  ${ n }_{ i }=\left\lfloor n{ \eta  }^{ -i } \right\rfloor $,   ${ r }_{ i }=r{ \eta  }^{ i }$;\\
            7 \quad \quad \quad \quad \quad  $\L \left( { \Lambda  }_{ s } \right) =\left\{ run\_ then\_ return\_ val\_ loss(\lambda ,\quad { r }_{ i })|\lambda \in { \Lambda  }_{ s } \right\} $;\\
            8 \quad \quad \quad \quad \quad ${ \Lambda  }_{ s }=top\_ k\left( { \Lambda  }_{ s },\quad \L \left( { \Lambda  }_{ s } \right) ,\quad \left\lfloor \frac { { n }_{ i } }{ \eta  }  \right\rfloor  \right)$ ;\\
            9 \quad \quad \quad \quad \textbf {end};\\
            10 \quad \textbf {end};\\
            11 \quad \textbf {output:} configuration $\lambda$ with lowest validation loss seen so far;\\
            \bottomrule 
        \end{tabular}
    }
    \label{Hyperband} 
\end{table}

\subsection{Genetic Algorithm}
\begin{figure}[H]
	\captionsetup{width=0.75\columnwidth}
	\centering\includegraphics[height=12cm,width=12.0cm]{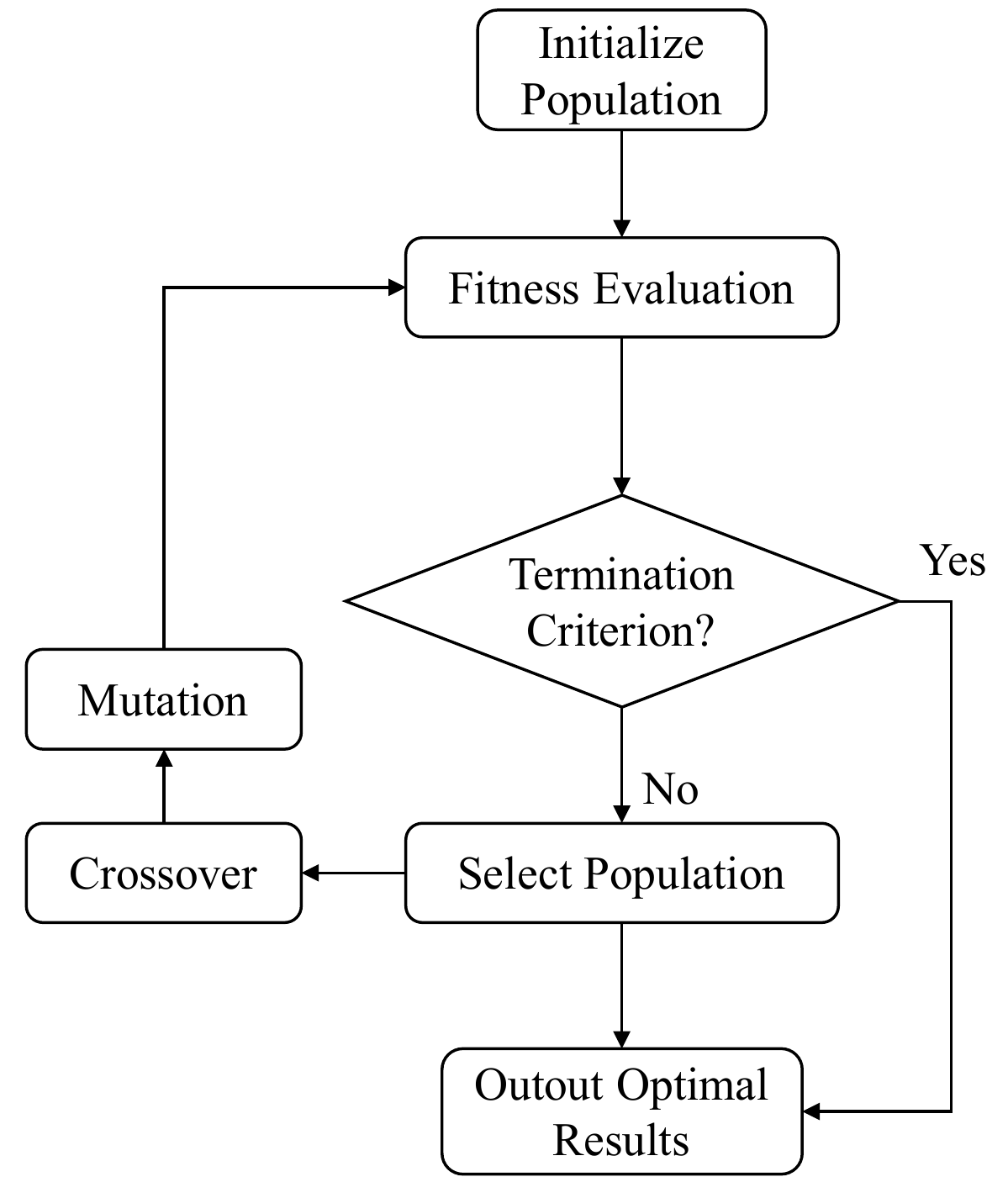} 
	\caption{ Genetic Algorithm}
	\label{GA}
\end{figure}

Based on the genetic and evolutionary processes of nature the GA is proposed. With the fitness function, the GA selects the parameters by selection, crossover, and mutation. In the process, the parameters with a good solution are retained and the parameters with a bad solution are eliminated, then the new populations are inherited. The solution of the new generation is better than the solution of the previous generation. We repeat the process until the optimization criteria are satisfied, in the end we can find the optimal solution. The process of the Genetic Algorithm is shown in Fig.~\ref{GA}.

The Genetic Algorithm has some advantages, compared to the traditional methods it's more efficient, it has good parallel capabilities, it can optimize both continuous and discrete functions and also multi-objective problems, and it is suitable for the problem with the large search space and  a large number of parameters. There are also some disadvantages for the Genetic Algorithm, The calculation for some problems is expensive due to the fitness value is calculated repeatedly, and the calculation time is relatively long. The Genetic Algorithm is hard to guarantee the optimality or the quality of the solution.

\subsection{Jaya Algorithm}
Recently the Jaya algorithm becomes a popular optimization method,  it can solve different types of optimization problems, for example, Rao used the Jaya algorithm to solve the constrained and unconstrained problems. The word Jaya comes from Sanskrit, which means victory.  When a particular problem is solved using the Jaya algorithm, we can get the best result,  and avoid the worst result. The Jaya algorithm is based on the population-based method, in the process, it is used to modify individual solutions repeatedly, and then use the best individual solutions to modify the population solution,  and in the end the best solution can be obtained. The Jaya algorithm is a gradient-free optimization algorithm. Compared with other algorithms, the Jaya algorithm isn't limited by hyper-parameters, and it only has two common parameters, that is, population size and the number of iterations. So Jaya algorithm is simple to find out the optimal solution and does not need to adjust any algorithm-specific parameters.
The Jaya algorithm has been applied for many fields.  such as B Chattopadhyay used the Jaya algorithm in the area of modern machining processes\cite{chattopadhyay1996line}, W Warid ad H Hizam used the Jaya algorithm to solve an optimal power flow solution \cite{warid2016optimal}. The Jaya algorithm was used for the heat exchangers by Rao and R Venkata \cite{rao2018multi}. 
 
The Jaya algorithm obtain the optimal values by minimizing the objective function,  we suppose the  objective function is defined, and it has n-dimensional factors, ${ v }_{ i }$ is the estimation value, i is the position of the candidate solution. The variables are updated using Eq.~(\ref{Jaya equation}).
\begin{equation}
\label{Jaya equation}  
  { v }_{ i+1 }={ v }_{ i }+{ r }_{ 1 }({ v }_{ b }-\left| { v }_{ i } \right| )-{ r }_{ 2 }({ v }_{ w }-\left| { v }_{ i } \right| )
\end{equation}  
where ${ v }_{ b }$ and ${ v }_{ w }$ are the best and worst solutions in the current population. ${ r }_{ 1 }$ and ${ r }_{ 2 }$ are  random numbers in the range of [0,1], which are used as scaling factors. The scaling factors attract the best solution and improve the worst of the update in each iteration. In this entire procedure,  the solution moves closer to the best result and moves away from the worst solution.
 
 With the Jaya algorithm, the objective function value gradually approaches the optimal solution by updating the value of the variable. In the process, the fitness of each candidate solution in the population is improved. The process of the Jaya algorithm is presented in following flowchart.
 \begin{figure}[H]
	\captionsetup{width=0.9\columnwidth}
	\centering\includegraphics[height=10cm,width=16.0cm]{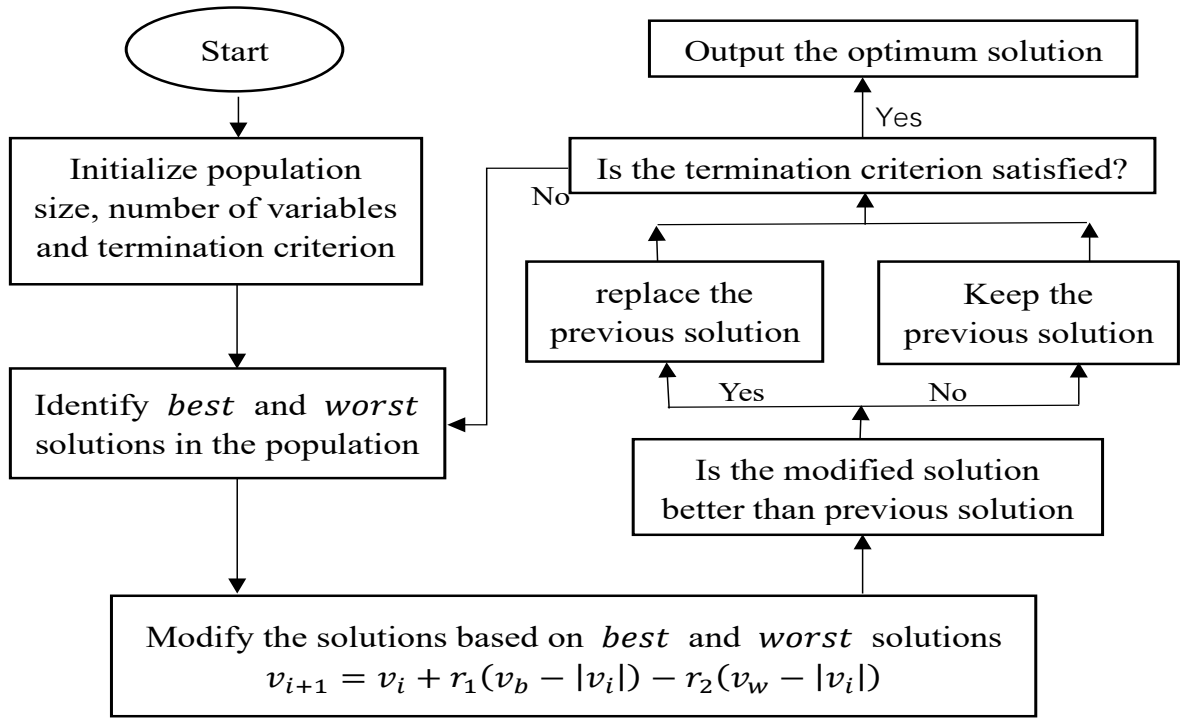}   
	\caption{ Flowchart of Jaya algorithm}
	\label{Jaya}
\end{figure}
\quad The performance of Jaya algorithm is reflected by the minimum function estimation.

\section{ The hybrid Jaya-Runge-Kutta Algorithm}
As is mentioned above, the shooting method can be sensitive to initial boundary conditions or unknown parameters in equations \cite{butcher2016numerical}, especially the guess of unknown initial values $\alpha$ and $\eta_{\infty}$ in Equations from Equations \ref{newfs} to \ref{newfsbd2}. The advancement of heuristic algorithms, however, opens a new door to deal with this type of problem easily and accurately.

\subsection{The Runge–Kutta method}
\begin{figure}[H]
	\captionsetup{width=0.9\columnwidth}
	\centering\includegraphics[height=8cm,width=14.0cm]{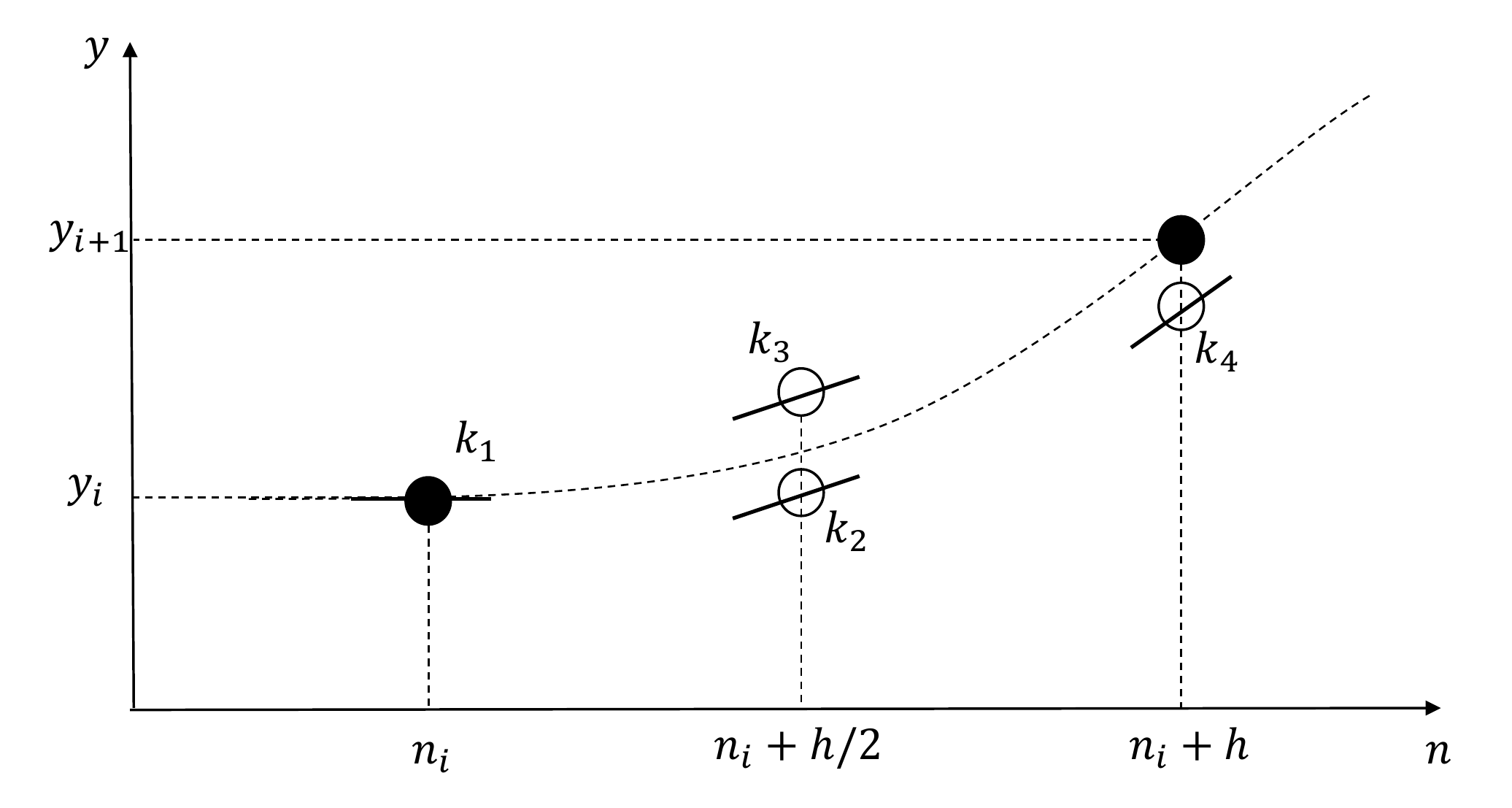} 
	\caption{ Fourth-order Runge-Kutta method}
	\label{runge-kutta}
\end{figure}

The Runge–Kutta method is a popular numerical analysis method, which solves the differential equation using a one-step nonlinear approximation. With its simplicity and efficiency, the method is widely used for solving the initial-value problems of differential equations. The Fourth-order Runge-Kutta method is shown in Fig.\ref{runge-kutta}, which is chosen as illustration. In each step, we evaluate the derivative four times, the first time is at the initial point, the second and third times are at trial midpoints, and the last time is at a trial endpoint , then the final function value is calculated with these derivatives, the final value is shown as a filled dot.
\quad Compared with other methods, the Fourth-order Runge-Kutta method is simple and robust scheme.

Consider the initial-boundary value problem from Equations \ref{newfs} to \ref{newfsbd2}, with a step-size h, the following recurrence formula can be obtained: 
\begin{equation}
\begin{aligned}
&{ k }_{ 1 }=hf({ \xi }_{ i },{ y }_{ i }),\\
& { k }_{ 2 }=hf({ \xi}_{ i }+\frac { h }{ 2 } ,{ y }_{ i }+\frac { { k }_{ 1 } }{ 2 } ),\\
& { k }_{ 3 }=hf({ \xi }_{ i }+\frac { h }{ 2 } ,{ y }_{ i }+\frac { { k }_{ 2 } }{ 2 } ),\\ 
&{ k }_{ 4 }=hf({ \xi }_{ i }+h,{ y }_{ i }+{ k }_{ 3 }),\\
& { y }_{ i+1 }={ y }_{ i }+\frac { 1 }{ 6 } ({ k }_{ 1 }+2{ k }_{ 2 }+2{ k }_{ 3 }+{ k }_{ 4 }),\\
& { \xi }_{ i+1 }={ \xi }_{ i }+h,
\end{aligned}
\end{equation}
where $f$ can be retrieved from Equation \ref{newfs}.

\subsection{The hybrid Jaya-Runge–Kutta method}
Unlike the classical shooting method, we start with a parameter interview as an initial guess to find appropriate values of $\alpha$ and $\eta_{\infty}$. We first begin by defining the objection function for Jaya algorithm.
\begin{equation}
min\quad F\left( \alpha ,{ \eta  }_{ \infty  } \right) = \parallel f_2\left( \alpha ,{ \eta  }_{ \infty  } \right) -f_2\left( { \eta  }_{ \infty  } \right);f_{ 3 }\left( \alpha ,{ \eta  }_{ \infty  } \right) -f_3\left( { \eta  }_{ \infty  }\right)   \parallel  
\end{equation}

 Since ${ RK }_{ 2 }\left( \alpha ,{ \eta  }_{ \infty  } \right)$ and ${ RK }_{ 3 }\left( \alpha ,{ \eta  }_{ \infty  } \right)$ are the numerical approximate boundary values, the objective function for Jaya algorithm can be reduced to:
 \begin{equation}\label{fitnessfun}
min\quad F\left( \alpha ,{ \eta  }_{ \infty  } \right) = \parallel { RK }_{ 2 }\left( \alpha ,{ \eta  }_{ \infty  } \right) -1;{ RK }_{ 3 }\left( \alpha ,{ \eta  }_{ \infty  } \right)   \parallel  
\end{equation}
Then the Jaya algorithm is deployed to minimize the objective function with 20 as
the size of initial population and a 100 maximum numbers of iterations to algorithm as termination criterion. 

$F$ is an objective function of $x=\left( \alpha ,{ \eta  }_{ \infty  } \right)$, supposing that $x_i,j$ is the estimation value of $j-$th valable for $i-$th competitor solution, where $j=1,2$. Therefore, $x_i=\left( \alpha_i ,{ \eta  }_{ \infty,i  } \right)$ is the position of $i-$th candidate solution. Let the best competitor solution $x_{best}=\left( \alpha_{best} ,{ \eta  }_{ \infty,best  } \right)$ gives the best estimation of $F$ in the present populace and the worst candidate solution $x_{worst}=\left( \alpha_{worst} ,{ \eta  }_{ \infty,worst } \right)$ obtains the worst value of $F$ in the present populace. Then the solution is modified based on the best and worst solution as:
 \begin{equation}\label{Jaya_rk}
x^{'}_{i,j}=x_{i,j}+Rand_1(x_{best}-\left | x_{i,j} \right |)-Rand_2(x_{worst}-\left | x_{i,j} \right |),
\end{equation} 
where $x^{'}_{i,j}$ is the update value of $x_{i,j}$. With Jaya optimization technique, the achieved solution moves closer to the finest result and starts moving away from the worst solution.

In what follows, the pseudo-code of the hybrid Jaya-Runge-Kutta algorithm is illustrated in Table.~\ref{Pseudo}, the method is used to solve the proposed optimization problem.

\begin{table}[H]
	\captionsetup{width=0.85\columnwidth}
	\caption{\textbf {Pseudo-code}} 
	\vspace{-0.3cm}
	\centering
    \resizebox{1\columnwidth}{!}{%
        \begin{tabular}{l} 
            \toprule 
            \textbf {Initialize}\\ 
            \quad randomly initialize $x_i=\left( { \alpha  }_{ i },{ \eta  }_{ \infty i } \right) $, i=1,2,...,N(Population size)\\
            \quad \textbf {while} Termination criterion is not met \textbf {do}\\
            \quad \quad calculate $R{ K }_{ 2 }({ \alpha  }_{ i },{ \eta  }_{ \infty i })$ and $R{ K }_{ 3 }({ \alpha  }_{ i },{ \eta  }_{ \infty i })$ by Runge-Kutta method\\
            \quad \quad evaluate the fitness function in Equation~\ref{fitnessfun}, i=1,2,...,N\\
            \quad \quad \quad \textbf {for} i in range(N) \textbf {do}\\
            \quad \quad \quad \quad find ${ x }_{ best }$ and ${ x }_{ worst }$.\\
            \quad \quad \quad \quad \textbf {for} j in range(2) \textbf {do}\\
            \quad \quad \quad \quad \quad update the variables by Eq.\ref{Jaya_rk}\\
            \quad \quad \quad \quad \quad next j(dimension of variables)\\
            \quad \quad \quad \quad \textbf {end}\\
            \quad \quad \quad \quad next i\\
            \quad \quad \quad \textbf {end}\\
            \quad \quad next generation until termination criterion satisfied\\
            \quad \textbf {end while}\\
            \bottomrule 
        \end{tabular}
    }
    \label{Pseudo} 
\end{table}

\section{Numerical results}
Next, the hybrid method is used to solve the transformed optimization problem, the Jaya algorithm is used to find the global minimum of the fitness function. For the best parameters $ \alpha$ and $ { \eta  }_{ \infty i }$ identified by Jaya  algorithm, Runge-Kutta method is then applied to solve the coupled differential equations. In order to test the ability of Jaya algorithm, the parameters identified by optimization methods are compared with results from open literatures \cite{zhang2009iterative,asaithambi2005solution}. The numerical results for the stream function $f$, velocity profile $f^{'}$, and skin friction coefficient $f^{''}$ are illustrated for specific flow problems. The velocity profile along the coordinates are compared with reference solutions from \cite{ahmad2017stochastic}, to show the accuracy of present method compared with other robust optimization methods.

 For different coefficients $\beta$ and ${ \beta  }_{ 0 }$, the Falkner-Skan problem can be divided into the following four categories:
\begin{itemize}
\item Blasius equation:${ \beta  }_{ 0 }=0.5, \beta =0$
\item Homann problem:${ \beta  }_{ 0 }=2,\beta =1$  
\item Accelerating flows:${ \beta  }_{ 0 }=1, \beta \geqslant 0$
\item Decelerating flows:${ \beta  }_{ 0 }=1,\beta <0$  
\end{itemize}
 
\subsection{Case1: Blasius equation with $\beta_{0}$=0.5 and $\beta$=0}
\begin{figure}[H]
	\captionsetup{width=0.9\columnwidth}
	\centering\includegraphics[height=10cm,width=14.0cm]{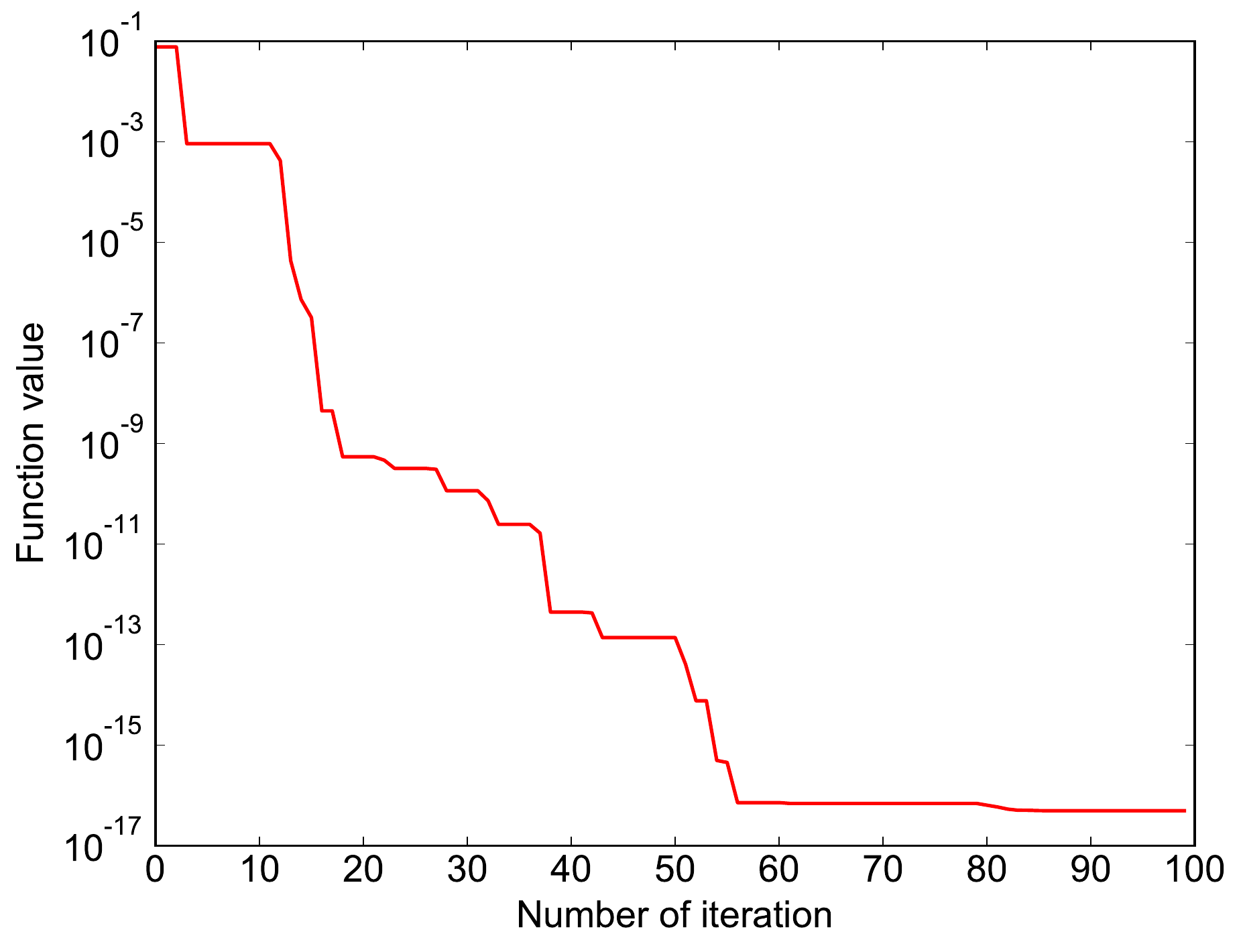}   
	\caption{ Convergence of fitness for ${ \beta  }_{ 0 }=0.5$ and ${ \beta  }$ = 0}
	\label{loss1}
\end{figure}

\begin{figure}[H]
	\captionsetup{width=0.9\columnwidth}
	\centering\includegraphics[height=9cm,width=13.0cm]{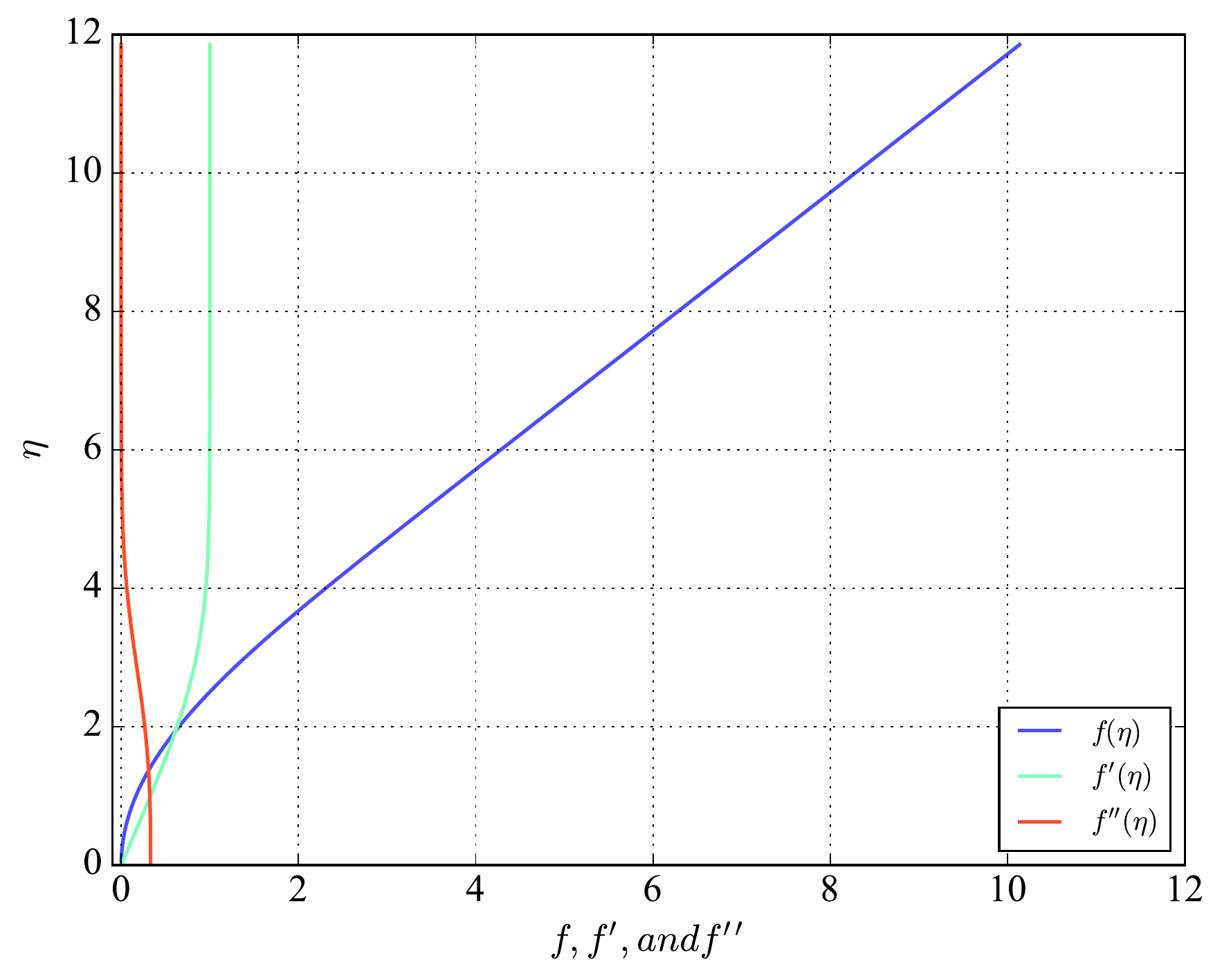}   
	\caption{  The stream function of Blasius equation and its derivatives
		corresponding to $\beta_{0}$=0.5, $\beta$=0}
	\label{case1:Blasius equation}
\end{figure}

\begin{figure}[H]
	\captionsetup{width=0.9\columnwidth}
	\centering\includegraphics[height=9cm,width=13.0cm]{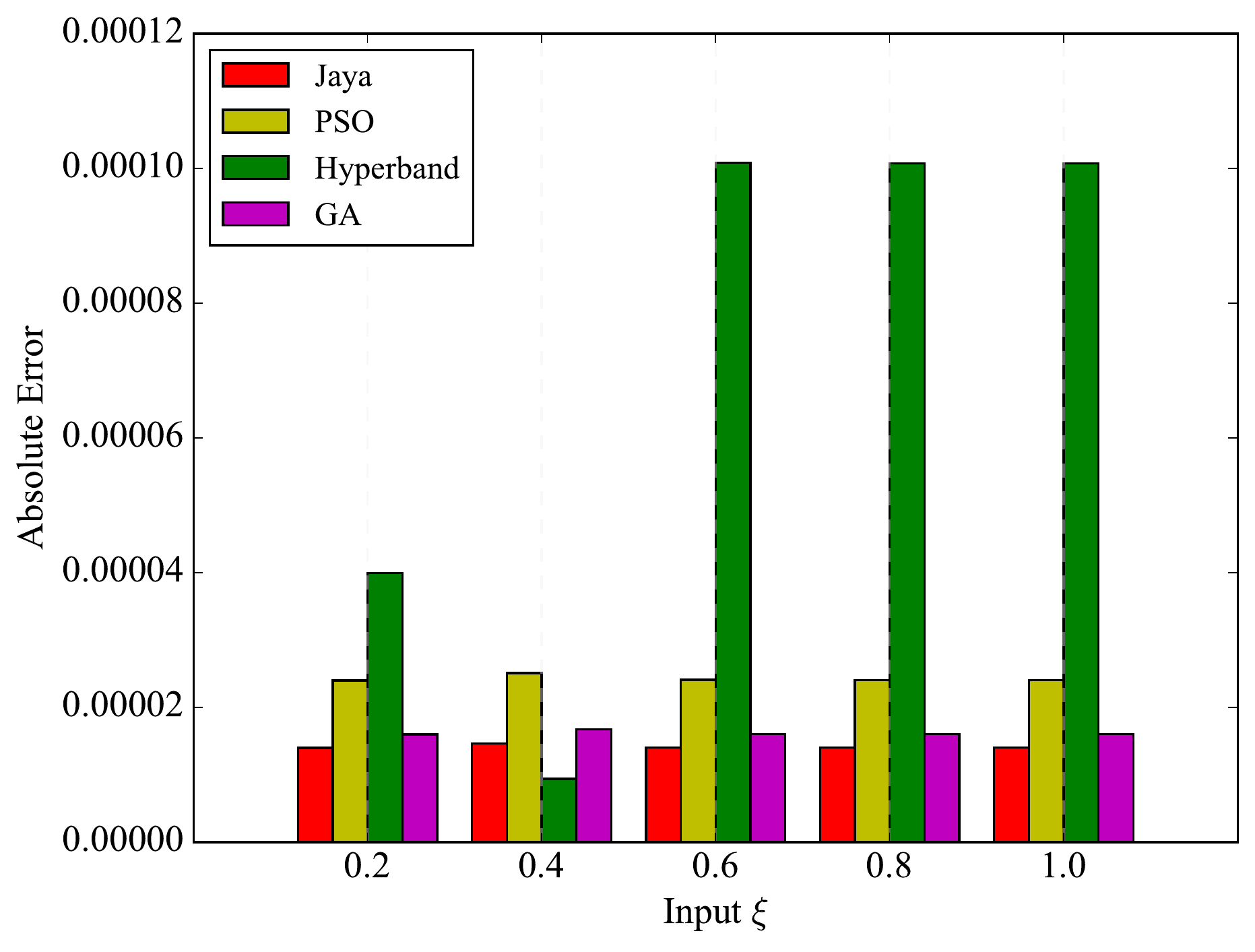}   
	\caption{Comparison of velocity of proposed results for Blasius equation}
	\label{compare1}
\end{figure}

In this case, Falkner–Skan equation is known as Blasius equation. The solution is obtained using the Runge-Kutta method integrated with the Jaya algorithm, the PSO algoritm, the Hyperband algorithm, the GA and reference solution from \cite{ahmad2017stochastic}. First, the convergence history graph of the fitness function over 100 iterations is presented in Figure \ref{loss1}. It is clearly the fitness function can reach its global minimum after some iterations, which proves the ability of Jaya algorithm as a global optimizer. 
Fig.~\ref{case1:Blasius equation}
shows the graph of stream function $f(\eta )$, its velocity and skin friction coefficient, where the velocity profile goes asymptotically
to 1 and the skin friction coefficient goes asymptotically
to 0, which are exactly the cases for Equations \ref{asymp} and \ref{bd3}. From the Table.~\ref{Table3}, the solutions from different optimization methods including Jaya, PSO, Hyperband and GA are compared with the reference solution. It is clearly that our method agrees very well with the reference solution. To be more clear, the absolution error for those listed optimization methods along the $\xi$ are presented Fig.~\ref{compare1}. It can be observed that the Jaya algorithm obtains the most stable and accurate results.

\begin{table}[H] 
	\captionsetup{width=0.85\columnwidth}
	\caption{Comparison of proposed results with reference solution for Blasius equation} 
	\vspace{-0.3cm}
	\centering 
	\resizebox{0.8\columnwidth}{!}{%
		\begin{tabular}{l c c c c c} 
			\toprule 
			\toprule 
			$\xi$&$f^{'}_{ref}$&Jaya&PSO&Hyperband&GA\\
			\midrule 
			0.1&0.3892354&0.38924351&0.38922149&0.38925129&0.38922613\\
			\midrule
			0.2&0.72239923&0.72241324&0.72237523&0.72235928&0.72238323\\
			\midrule
			0.3&0.91901667&0.91903214&0.91899015&0.91903766&0.91899899\\
			\midrule
			0.4&0.98646615&0.9864808&0.98644101&0.98647555&0.98644939\\
			\midrule
			0.5&0.99877239&0.99878654&0.99874812&0.99877457&0.99875621\\
			\midrule
			0.6&0.9999285&0.99994256&0.9999044&1.00002937&0.99991243\\
			\midrule
			0.7&0.99998388&0.99999793&0.99995978&1.00008464&0.99996781\\
			\midrule
			0.8&0.99998522&0.99999928&0.99996113&1.00008598&0.99996916\\
			\midrule
			0.9&0.99998524&0.99999929&0.99996115&1.000086&0.99996918\\
			\midrule
			1&0.99998524&0.99999929&0.99996115&1.000086&0.99996918\\
			\bottomrule 
		\end{tabular}
	}
	\label{Table3} 
\end{table}

\subsection{Case2: Homann flow with ${ \beta  }_{ 0 }=2,\beta =1$}
\begin{figure}[H]
	\captionsetup{width=0.9\columnwidth}
	\centering\includegraphics[height=10cm,width=14.0cm]{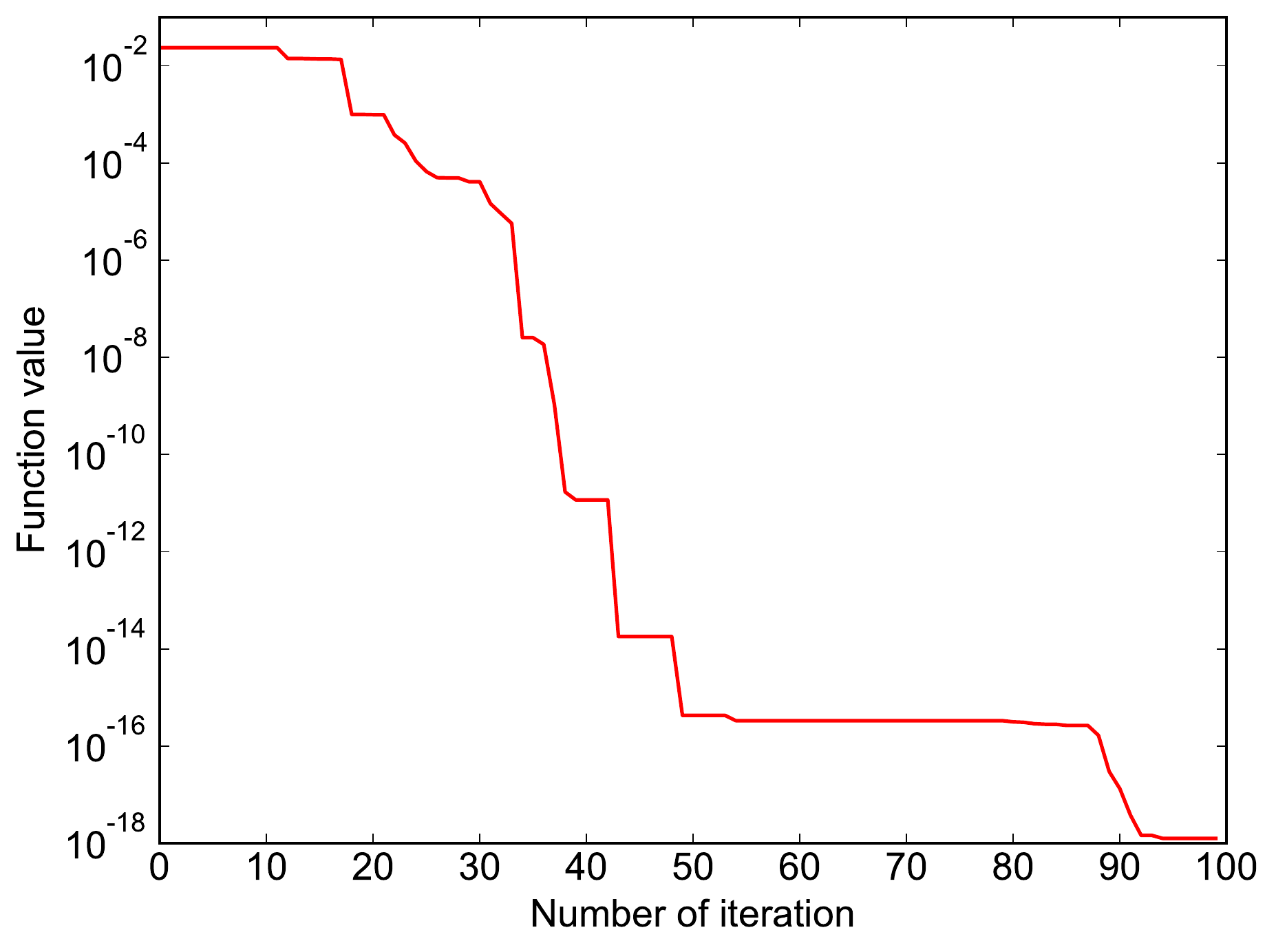}   
	\caption{ Convergence of fitness for ${ \beta  }_{ 0 }=2$ and ${ \beta  }$ = 1}
	\label{loss2}
\end{figure}

\begin{figure}[H]
	\captionsetup{width=0.9\columnwidth}
	\centering\includegraphics[height=9cm,width=13.0cm]{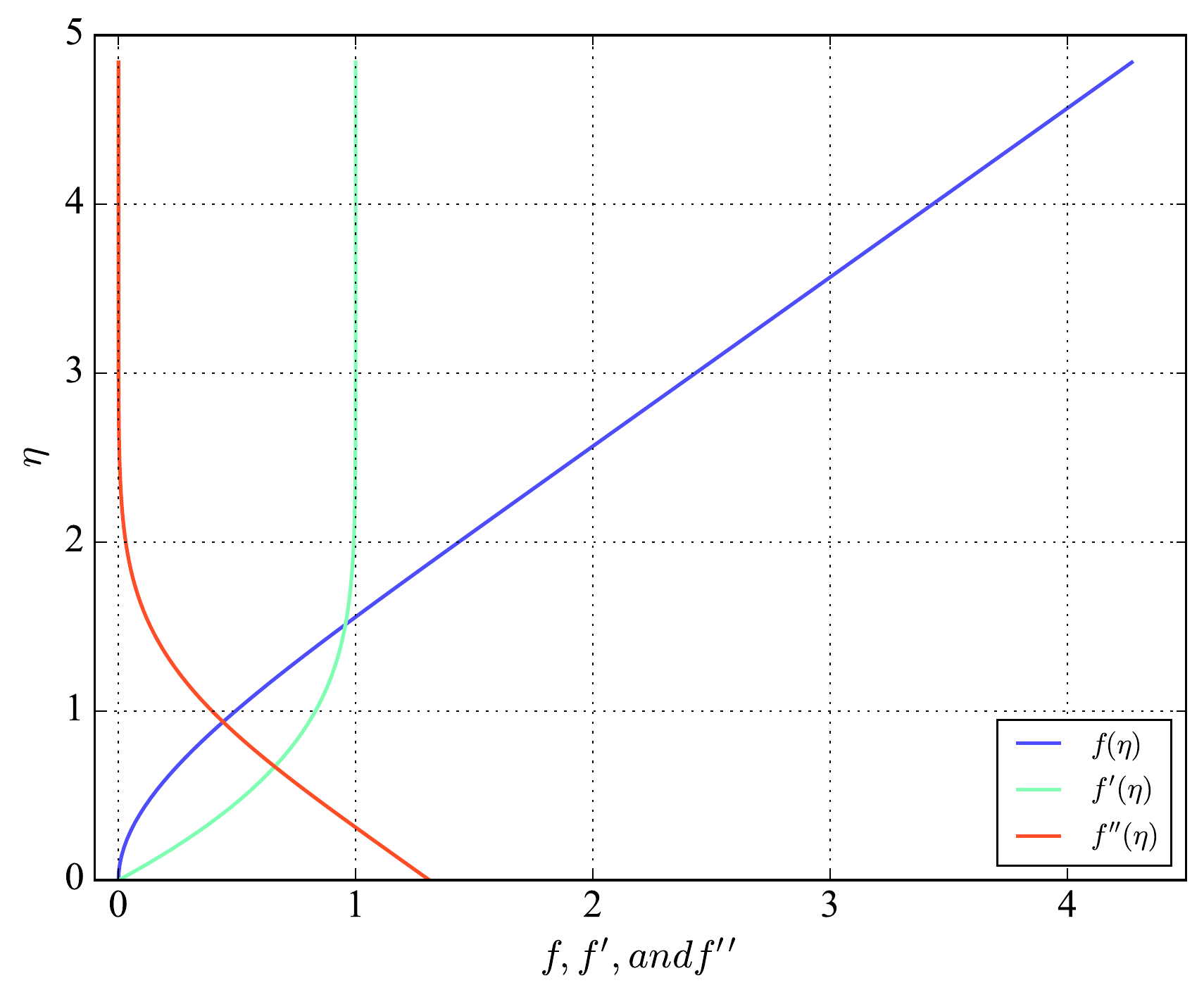}   
	\caption{  The stream function for Homann flow and its derivatives
		corresponding to $\beta_{0}$=2, $\beta$=1}
	\label{case4}
\end{figure}

In this section, the Falkner–Skan problem is considered as Homann steady flow. Also, the solution is obtained with the help of Jaya, PSO, Hyperband, GA integrated Runge-Kutta method. First, the convergence history graph of the fitness function over 100 iterations is presented in Figure \ref{loss2}. Also, the value of fitness function drops down significantly. 
Fig.~\ref{case4}
shows the graph of stream function $f(\eta )$, its velocity and skin friction coefficient gained with Jaya optimization method. From Table.~\ref{Table4}, the solutions from different optimization methods including Jaya, PSO, Hyperband and GA are compared with the reference solution \cite{ahmad2017stochastic}. Still our method obtains the most agreeable results. Moreover, the absolution errors for those listed optimization methods along the $\xi$ are presented Fig.~\ref{compare2}. It can be concluded that the Jaya algorithm obtains the most stable and accurate results and far outweighs the other optimization method.
 
\begin{table}[H] 
	\captionsetup{width=0.85\columnwidth}
	\caption{Comparison of proposed results with reference solution for Homann flow} 
	\vspace{-0.3cm}
	\centering 
	\resizebox{0.8\columnwidth}{!}{%
		\begin{tabular}{l c c c c c} 
			\toprule 
			\toprule 
			$\xi$&$f^{'}_{ref}$&Jaya&PSO&Hyperband&GA\\
			\midrule 
			0.1&0.51841208&0.51841111&0.5183927&0.51837284&0.51842661\\
			\midrule
			0.2&0.81685809&0.81685613&0.81681886&0.81677866&0.81688751\\
			\midrule
			0.3&0.9483164&0.94831334&0.94825528&0.94819263&0.94836223\\
			\midrule
			0.4&0.98969146&0.98968714&0.98960504&0.98951645&0.98975628\\
			\midrule
			0.5&0.99859956&0.99859383&0.998485&0.99836758&0.99868548\\
			\midrule
			0.6&0.99988001&0.99987281&0.99973603&0.99958845&0.999988\\
			\midrule
			0.7&1.00000237&0.99999368&0.99982862&0.99965052&1.00013269\\
			\midrule
			0.8&1.00001151&1.00000133&0.99980793&0.99979926&1.00016419\\
			\midrule
			0.9&1.00001353&1.00000186&0.99978012&0.99964087&1.00018858\\
			\midrule
			1&1.00001527&1.00000211&0.99975203&0.9996822&1.00021269\\

			\bottomrule 
		\end{tabular}
	}
	\label{Table6} 
\end{table}

\begin{figure}[H]
	\captionsetup{width=0.9\columnwidth}
	\centering\includegraphics[height=9cm,width=13.0cm]{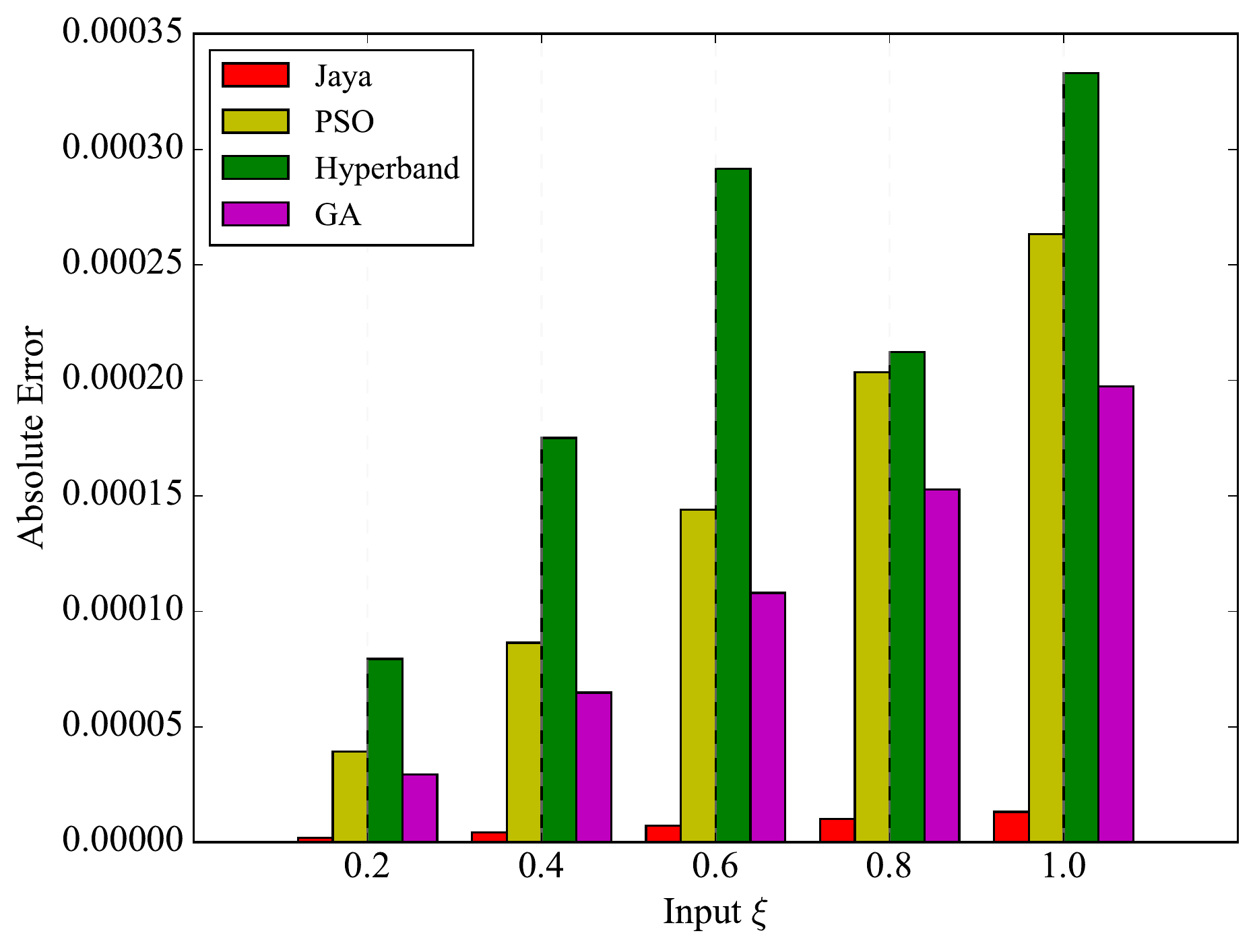}   
	\caption{ Comparison of velocity of proposed results for Homann flow}
	\label{compare4}
\end{figure}

\subsection{Case3: Accelerating flows with ${ \beta  }_{ 0 }=1, \beta \geqslant 0$}
\begin{figure}[H]
	\captionsetup{width=0.9\columnwidth}
	\centering\includegraphics[height=10cm,width=14.0cm]{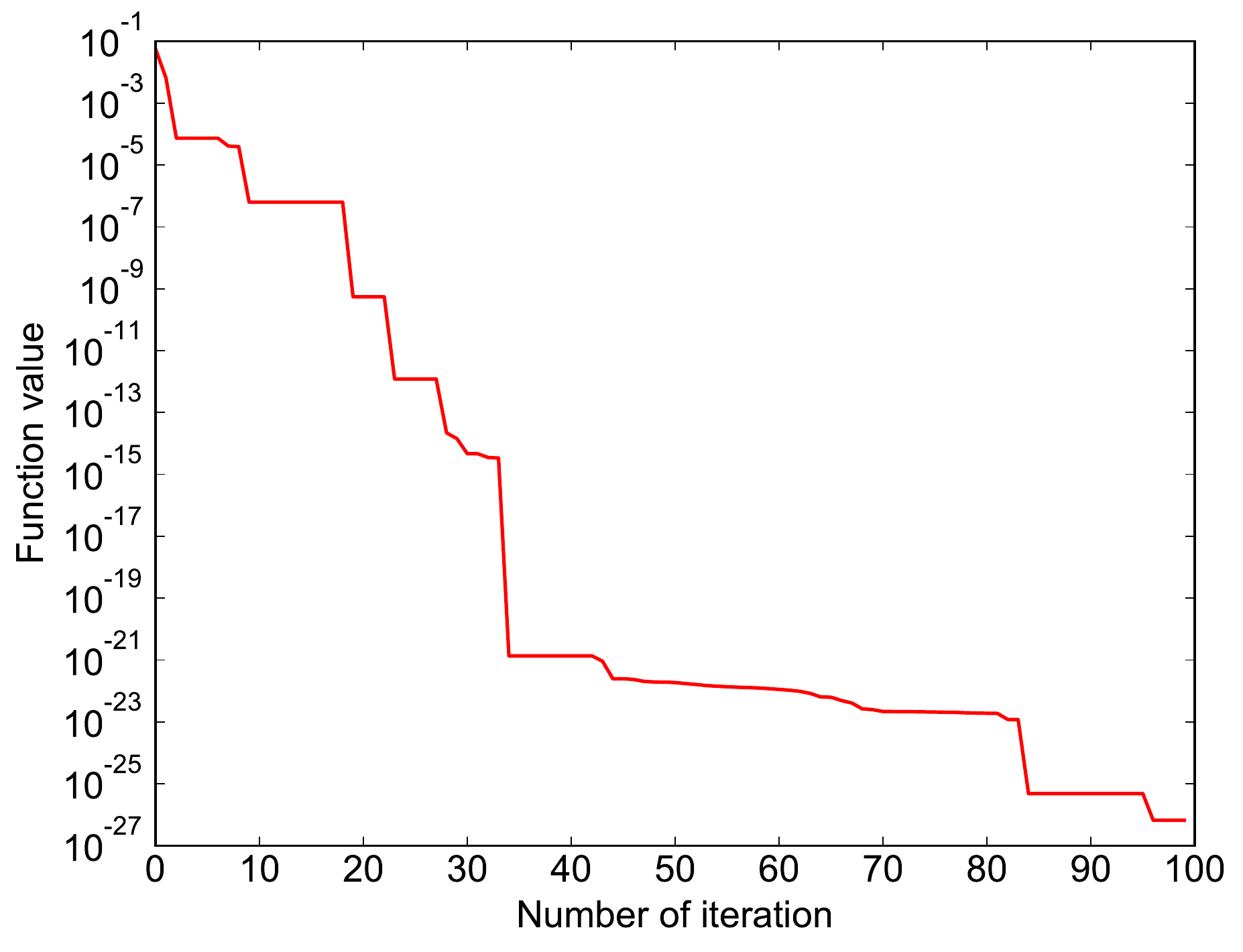}   
	\caption{ Convergence of fitness for ${ \beta  }_{ 0 }=1$ and ${ \beta  }$ = 0}
	\label{loss6}
\end{figure}

\begin{figure}[H]
	\captionsetup{width=0.9\columnwidth}
	\centering\includegraphics[height=10cm,width=14.0cm]{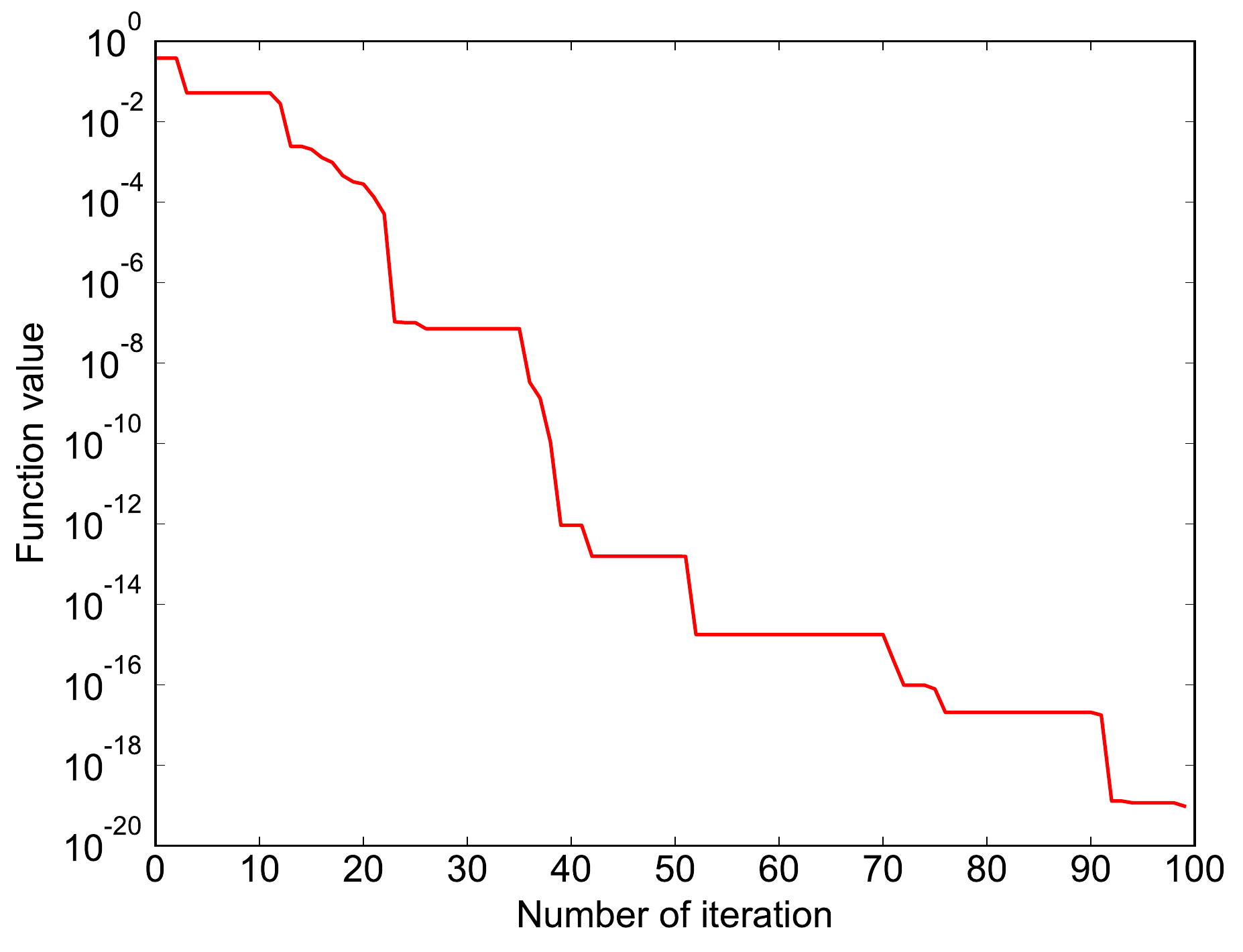}   
	\caption{ Convergence of fitness for ${ \beta  }_{ 0 }=1$ and ${ \beta  }$ = 0.5}
	\label{loss5}
\end{figure}

\begin{figure}[H]
	\captionsetup{width=0.9\columnwidth}
	\centering\includegraphics[height=10cm,width=14.0cm]{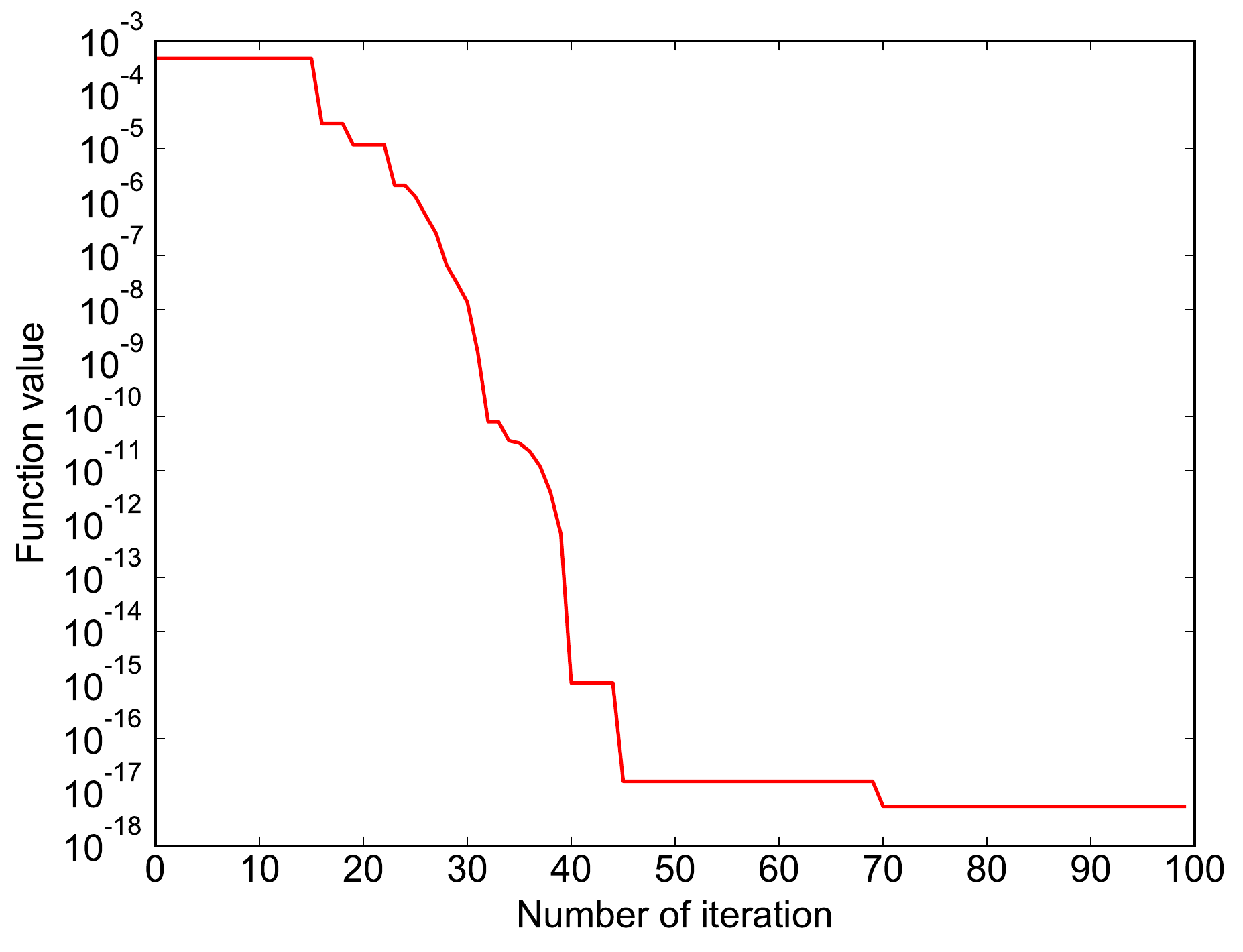}   
	\caption{ Convergence of fitness for ${ \beta  }_{ 0 }=1$ and ${ \beta  }$ =1 }
	\label{loss3}
\end{figure}

\begin{figure}[H]
	\captionsetup{width=0.9\columnwidth}
	\centering\includegraphics[height=10cm,width=14.0cm]{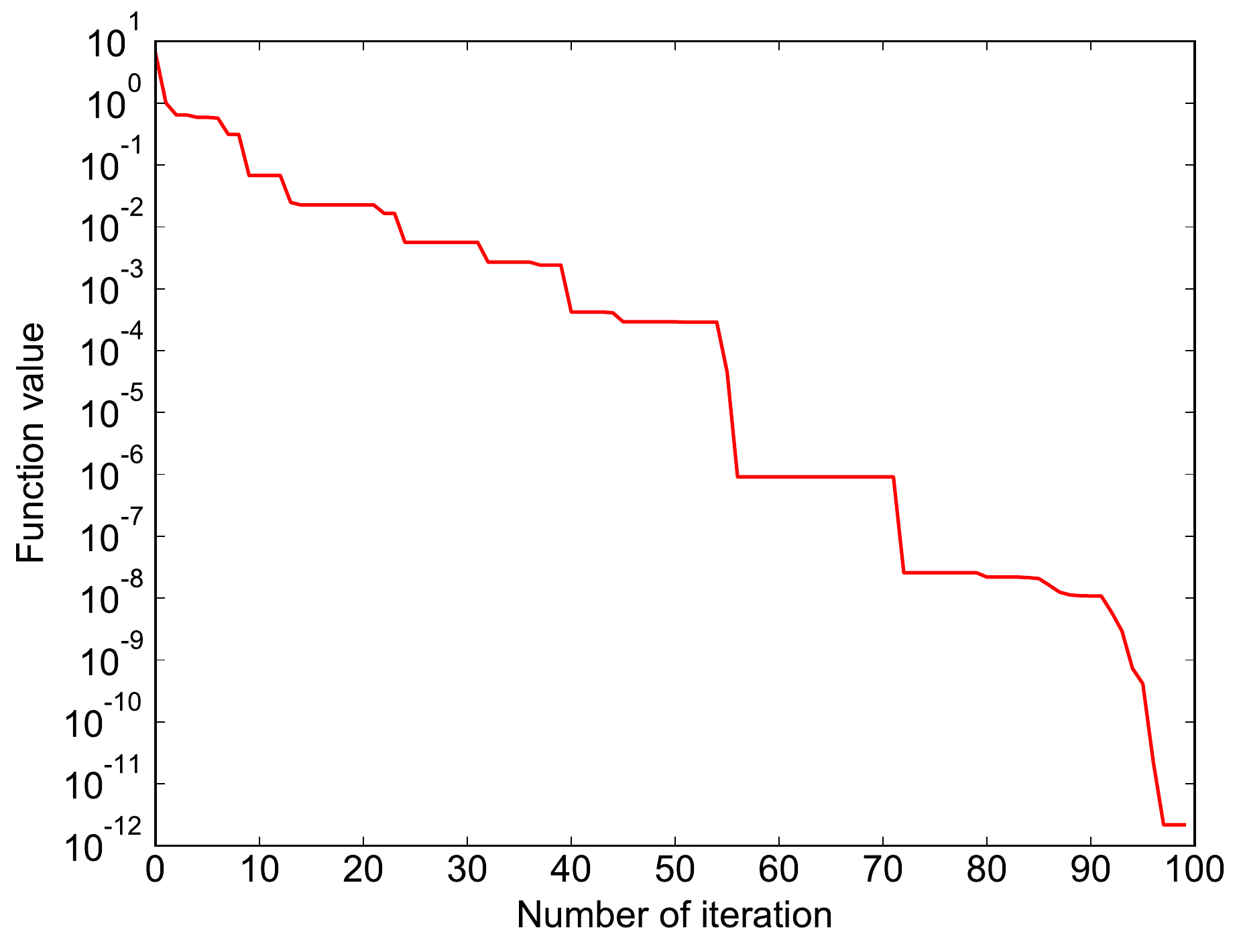}   
	\caption{ Convergence of fitness for ${ \beta  }_{ 0 }=1$ and ${ \beta  }$ =2}
	\label{loss4}
\end{figure}

The accelerating flows are studied in this section. Different $\beta$ values are selected. First, the convergence history graph of the fitness function over 100 iterations for different parameters is presented from Figures \ref{loss6} to \ref{loss4}. Also, the value of fitness function can reach a very low level at the end of the iterations.

\begin{figure}[H]
	\captionsetup{width=0.9\columnwidth}
	\centering\includegraphics[height=9cm,width=13.0cm]{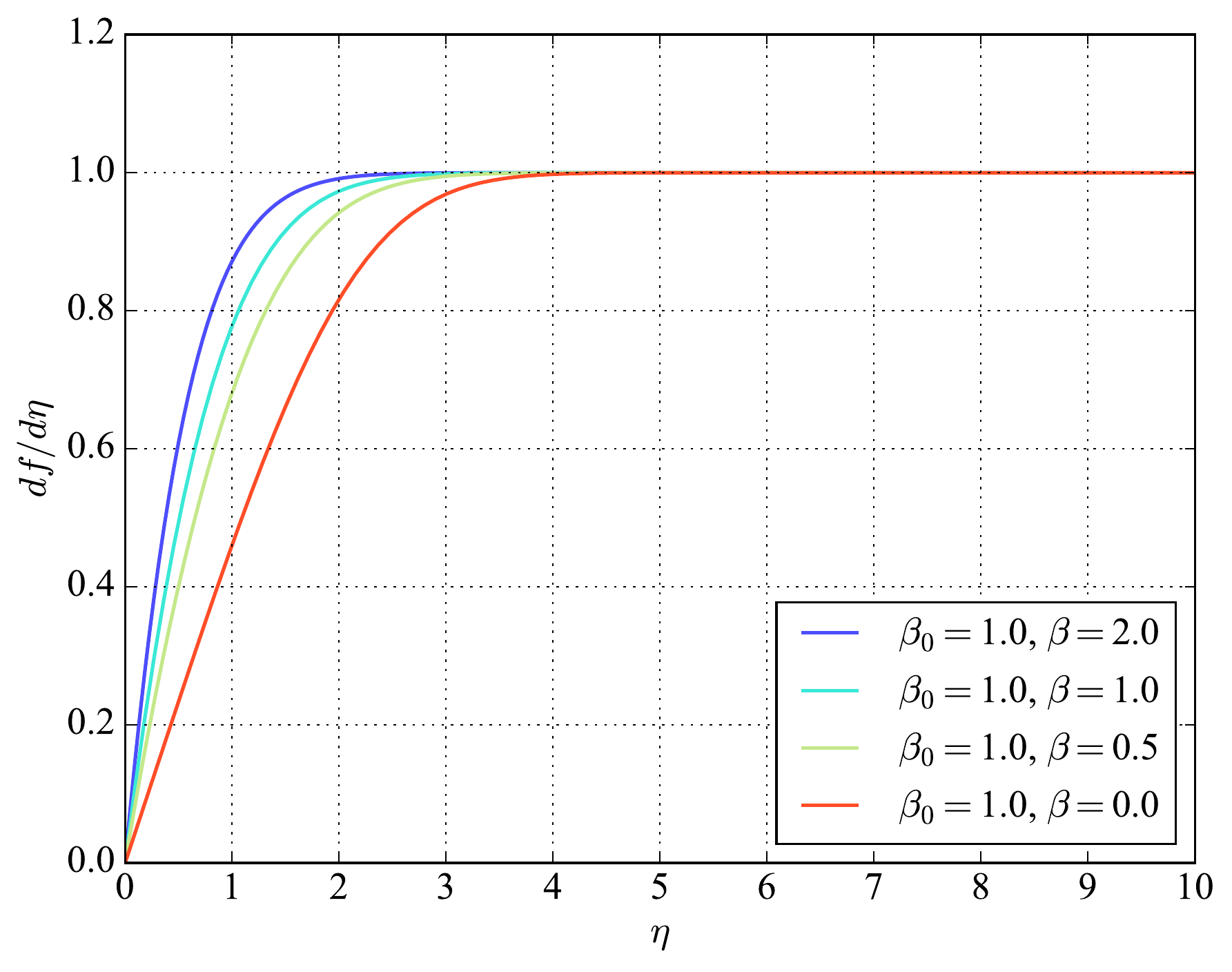}   
	\caption{ The velocity profile corresponding to different  $\beta \in \left[ 0,2 \right] $ when ${ \beta  }_{ 0 }=1$}
	\label{beta2}
\end{figure}

Fig.~\ref{beta2}
shows the velocity profile for those selected cases gained with Jaya optimization method. With the increasing of $\eta$, it is observed
 that the horizontal velocity profiles go asymptotically
to 1, which verifies the asymptotical boundary condition shown in Equation~\ref{asymp}. In addition, for acceleration flows, we can observe that as the
parameter $\beta$ increases, the boundary layer thickness
increases, and eventually tends to one as the distance increases
from the initial boundary. 

\subsubsection{ Hiemenz flow problem for $\beta_{0}$=1 and $\beta$=1}
\begin{figure}[H]
	\captionsetup{width=0.9\columnwidth}
	\centering\includegraphics[height=9cm,width=13.0cm]{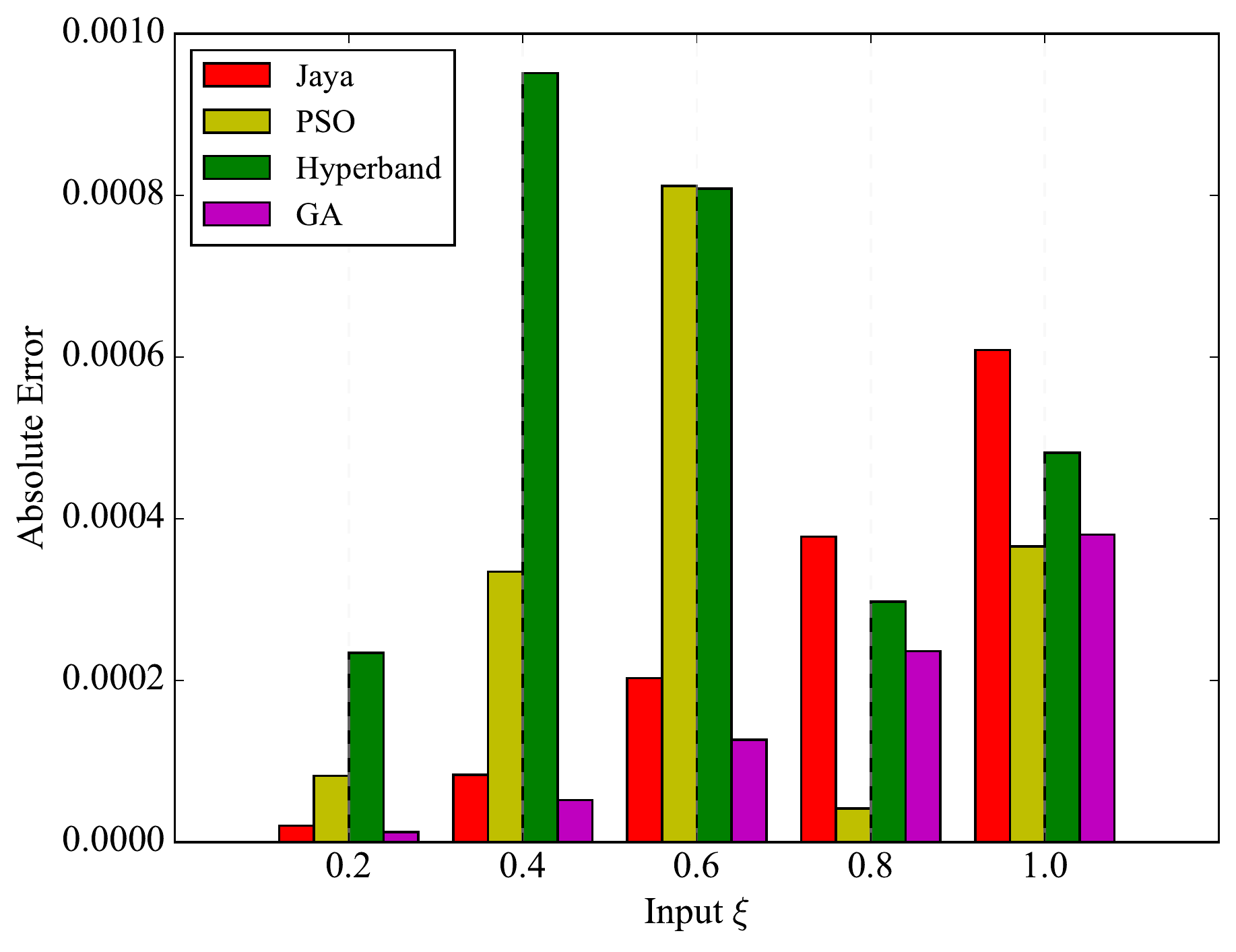}   
	\caption{Comparison of velocity of proposed results for Hiemenz flow}
	\label{compare2}
\end{figure}

To be more specific, two classic cases Hiemenz flow and Homann axisymmetric stagnation flow, are in comparison with Jaya, PSO, Hyperband, GA, and the reference solution \cite{ahmad2017stochastic}. In Tables~\ref{Table4}, the results obtained with Jaya method agrees pretty well with the reference solution \cite{ahmad2017stochastic}. To be more clear, the absolution error for those listed optimization methods along the $\xi$ are presented Fig.~\ref{compare2}. For this case, though Jaya optimizer is not the best, it still gives very accurate results at most points.
\begin{figure}[H]
	\captionsetup{width=0.9\columnwidth}
	\centering\includegraphics[height=9cm,width=13.0cm]{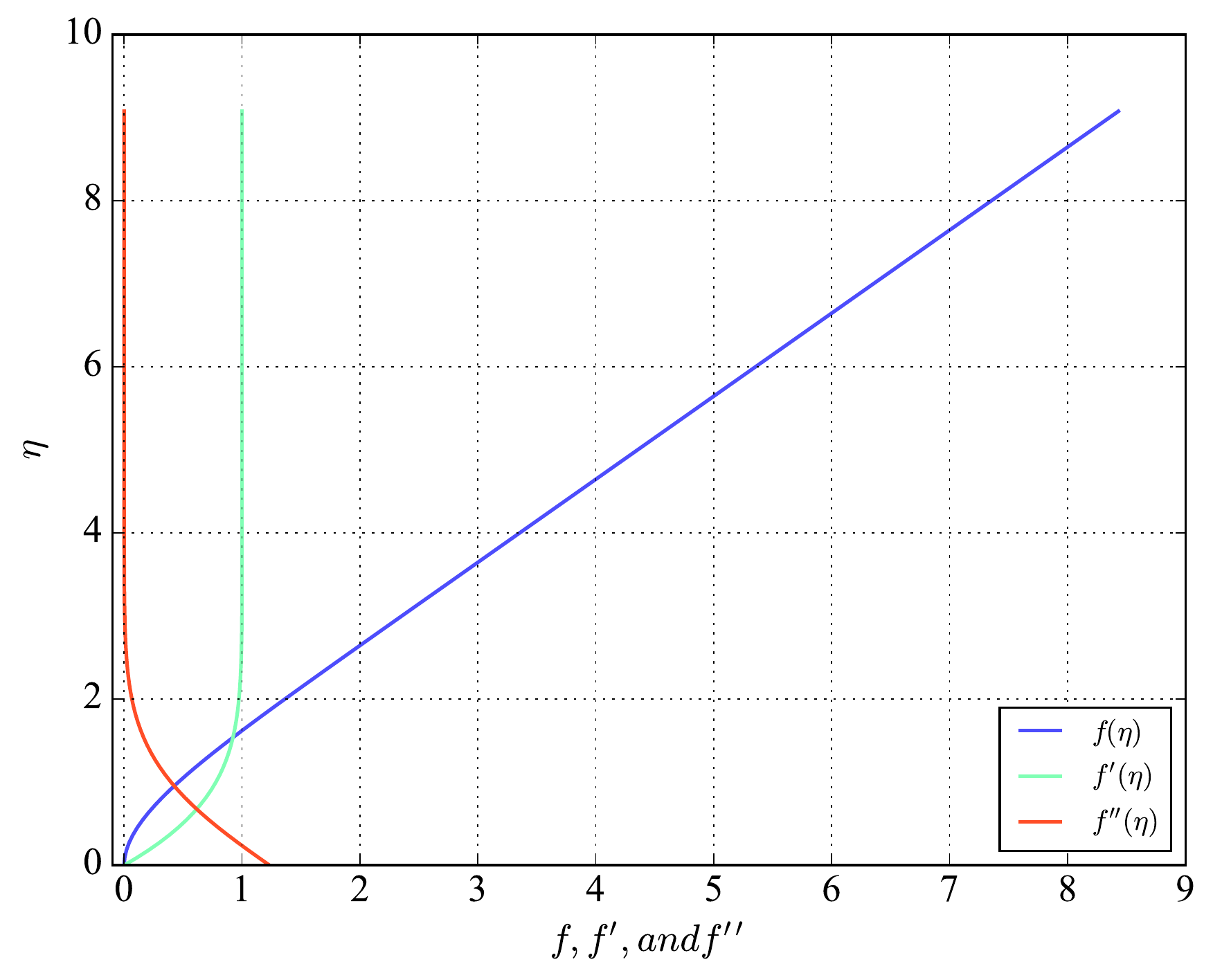}   
	\caption{  The stream function for Hiemenz flow and its derivatives
		corresponding to $\beta_{0}$=1, $\beta$=1}
	\label{case2}
\end{figure} 

\begin{table}[H] 
	\captionsetup{width=0.85\columnwidth}
	\caption{Comparison of proposed results with reference solution for Hiemenz flow} 
	\vspace{-0.3cm}
	\centering 
	\resizebox{0.8\columnwidth}{!}{%
		\begin{tabular}{l c c c c c} 
			\toprule 
			\toprule 
			$\xi$&$f^{'}_{ref}$&Jaya&PSO&Hyperband&GA\\
			\midrule 
			0.1&0.73864188&0.73864962&0.73861094&0.73855389&0.73863705\\
			\midrule
			0.2&0.9583553&0.95837589&0.95827293&0.95812107&0.95834243\\
			\midrule
			0.3&0.99623432&0.99627939&0.99605404&0.99572166&0.99620615\\
			\midrule
			0.4&0.99975182&0.99983548&0.99941718&0.99880024&0.99969953\\
			\midrule
			0.5&0.99987682&1.00001317&0.99933143&0.99892602&0.9997916\\
			\midrule
			0.6&0.99982269&1.00002568&0.9990108&0.99901422&0.99969582\\
			\midrule
			0.7&0.99975239&1.00003594&0.99911834&0.99912817&0.99957517\\
			\midrule
			0.8&0.99966988&1.00004792&0.99962806&0.99937203&0.99943362\\
			\midrule
			0.9&0.99957522&1.00006166&0.99973003&0.99944605&0.99927122\\
			\midrule
			1&0.9994684&1.00007717&0.99983429&0.99995043&0.99908797\\
			\bottomrule 
		\end{tabular}
	}
	\label{Table4} 
\end{table}

\subsubsection{Homann axisymmetric stagnation flow for $\beta_{0}$=1 and $\beta$=0.5}
\begin{figure}[H]
	\captionsetup{width=0.9\columnwidth}
	\centering\includegraphics[height=9cm,width=13.0cm]{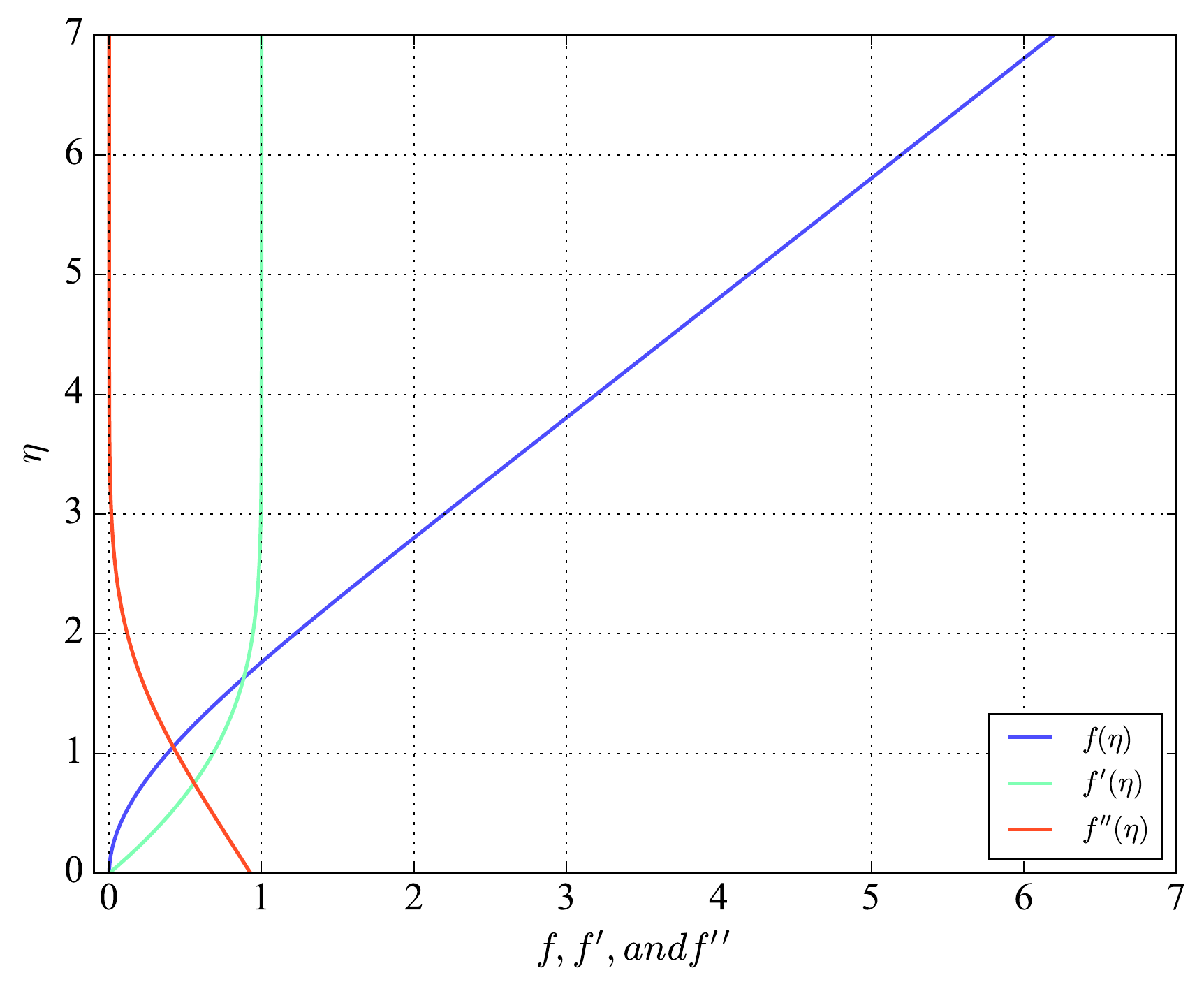}   
	\caption{ The stream function for Homann axisymmetric stagnation flow and its derivatives
		corresponding to $\beta_{0}$=1, $\beta$=0.5}
	\label{case3}
\end{figure}

\begin{figure}[H]
	\captionsetup{width=0.9\columnwidth}
	\centering\includegraphics[height=9cm,width=13.0cm]{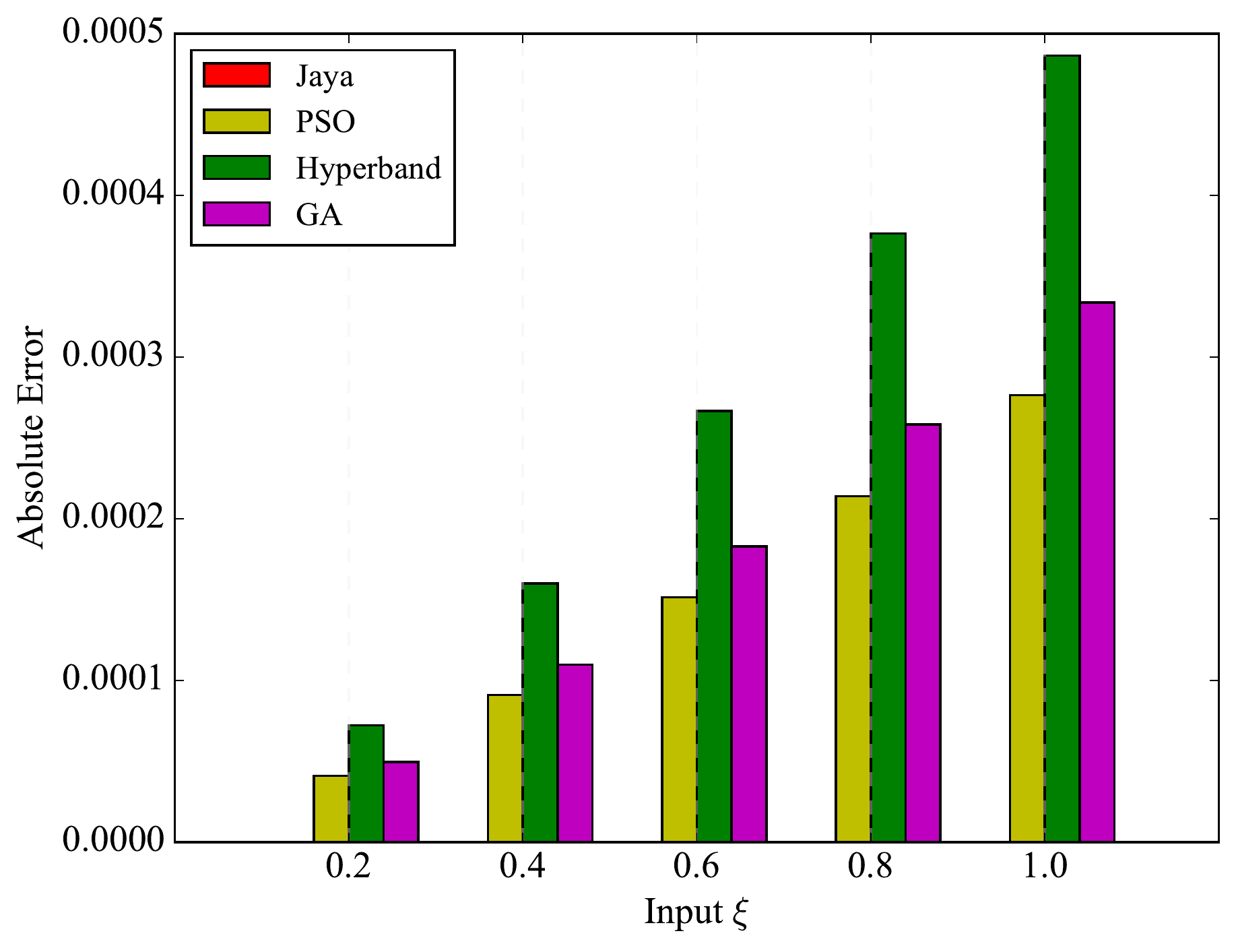}   
	\caption{Comparison of velocity of proposed results for Homann axisymmetric stagnation flow}
	\label{compare3}
\end{figure}

For the Homann axisymmetric stagnation flow, the stream function $f(\eta )$, its velocity and skin friction coefficient gained with Jaya optimization method are first shown in Figure \ref{case3}. Also, it can be observed that with the increase of $\eta$, the velocity profile goes asymptotically
to 1, which verifies the asymptotical boundary condition shown in Equation~\ref{asymp} and the skin friction coefficient goes asymptotically
to 0, which is just the case for Equation \ref{bd3}. In detail, the velocity profiles gained by all optimization methods are compared with reference solution \cite{ahmad2017stochastic} shown in Table~\ref{Table5}. The results obtained by Jaya method are in excellent agreement with the reference solution. For better comparison, the absolution error for those listed optimization methods along the $\xi$ is first shown in Fig \ref{compare3}. The hybrid Jaya Runge-Kutta method gains results with the absolute errors almost zero, it can hardly be seen from the error bar graph.

\begin{table}[H] 
	\captionsetup{width=0.85\columnwidth}
	\caption{Comparison of proposed results with reference solution for Homann axisymmetric stagnation flow} 
	\vspace{-0.3cm}
	\centering 
	\resizebox{0.8\columnwidth}{!}{%
		\begin{tabular}{l c c c c c} 
			\toprule 
			\toprule 
			$\xi$&$f^{'}_{ref}$&Jaya&PSO&Hyperband&GA\\
			\midrule 
			0.1&0.52720714&0.52720713&0.52718683&0.52717143&0.52718263\\
			\midrule
			0.2&0.82558612&0.82558615&0.82554498&0.82551377&0.82553647\\
			\midrule
			0.3&0.95299636&0.95299639&0.9529321&0.95288335&0.9529188\\
			\midrule
			0.4&0.99120424&0.99120421&0.9911132&0.99104412&0.99109436\\
			\midrule
			0.5&0.998896&0.998895&0.99877527&0.99868369&0.9987503\\
			\midrule
			0.6&0.99990906&0.99990904&0.99975742&0.99964238&0.99972605\\
			\midrule
			0.7&0.999995&0.999994&0.99981215&0.99967344&0.99977432\\
			\midrule
			0.8&0.99999962&0.9999996&0.99978552&0.9996231&0.99974123\\
			\midrule
			0.9&0.99999975&0.99999979&0.9997544&0.99956826&0.99970363\\
			\midrule
			1&0.99999972&0.99999979&0.99972311&0.99951327&0.99966588\\
			\bottomrule 
		\end{tabular}
	}
	\label{Table5} 
\end{table}

\subsection{Case4: Decelerating flows with ${ \beta  }_{ 0 }=1,\beta <0$}
\begin{figure}[H]
	\captionsetup{width=0.9\columnwidth}
	\centering\includegraphics[height=10cm,width=14.0cm]{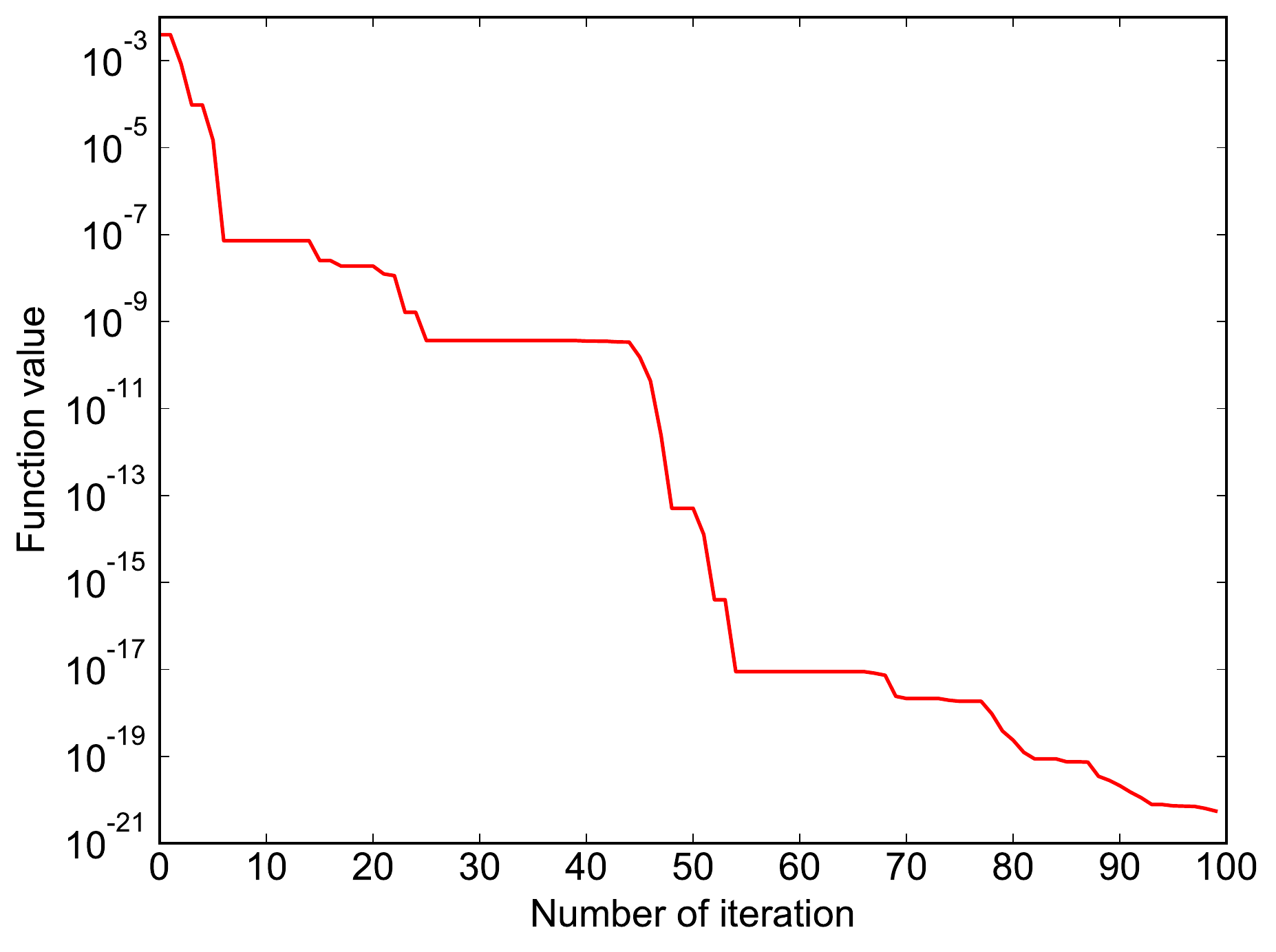}   
	\caption{ Convergence of fitness for ${ \beta  }_{ 0 }=1$ and ${ \beta  }$ = -0.1}
	\label{loss7}
\end{figure}

\begin{figure}[H]
	\captionsetup{width=0.9\columnwidth}
	\centering\includegraphics[height=10cm,width=14.0cm]{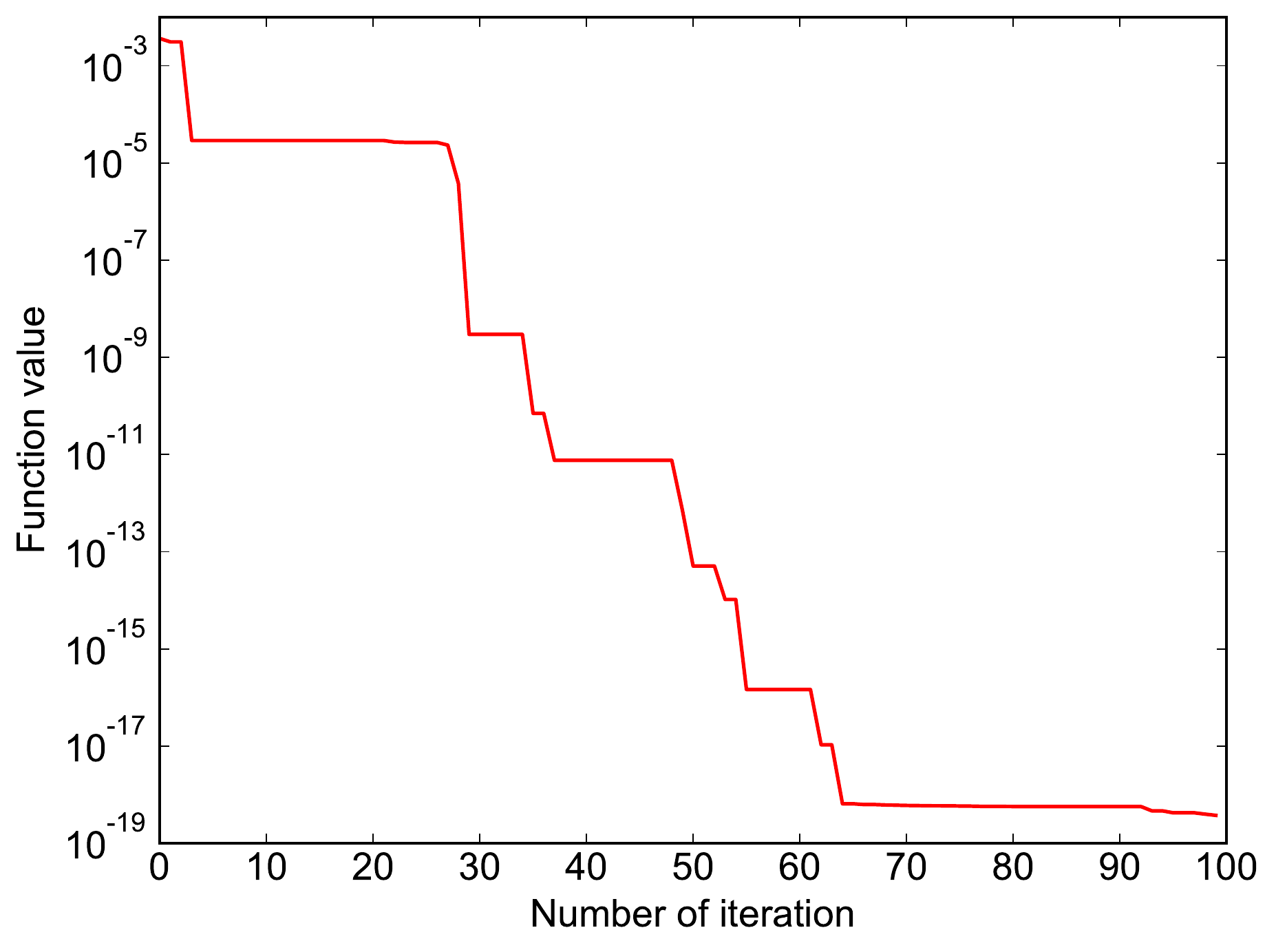}   
	\caption{ Convergence of fitness for ${ \beta  }_{ 0 }=1$ and ${ \beta  }$ = -0.15}
	\label{loss8}
\end{figure}

\begin{figure}[H]
	\captionsetup{width=0.9\columnwidth}
	\centering\includegraphics[height=10cm,width=14.0cm]{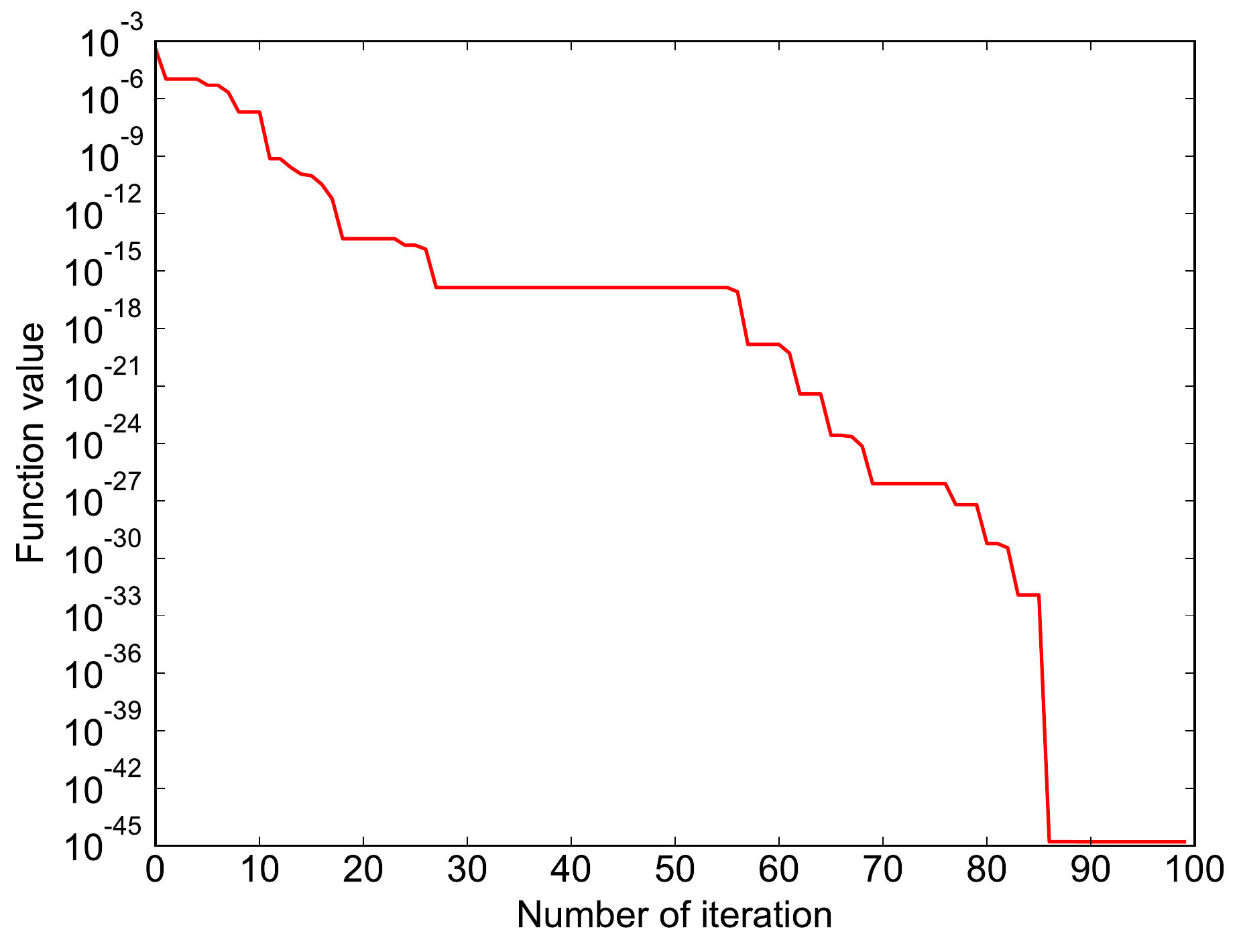}   
	\caption{ Convergence of fitness for ${ \beta  }_{ 0 }=1$ and ${ \beta  }$ = -0.18}
	\label{loss9}
\end{figure}

\begin{figure}[H]
	\captionsetup{width=0.9\columnwidth}
	\centering\includegraphics[height=10cm,width=14.0cm]{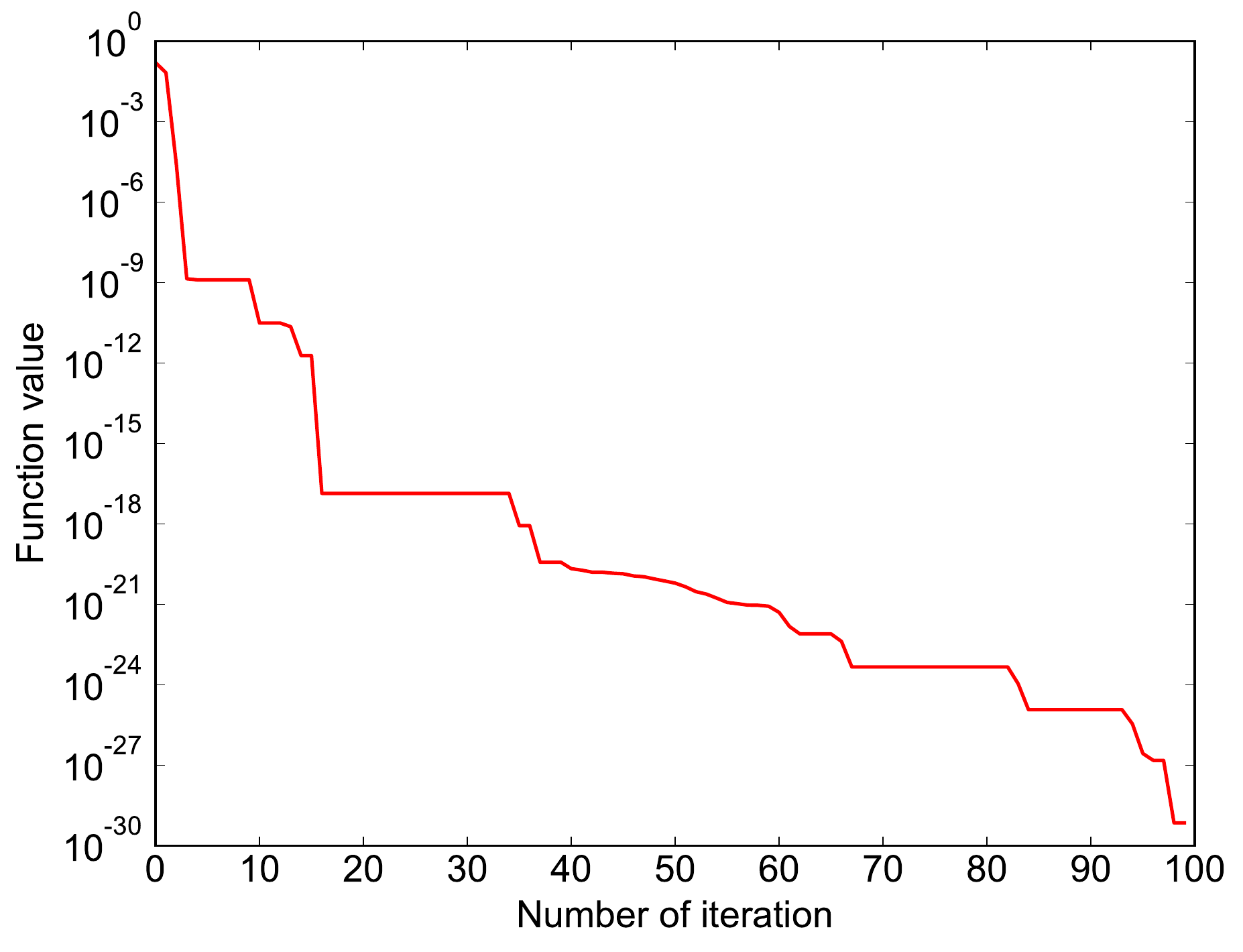}   
	\caption{ Convergence of fitness for ${ \beta  }_{ 0 }=1$ and ${ \beta  }$ =-0.1988}
	\label{loss10}
\end{figure}

Finally, the decelerating flows are studied in this section. Also various $\beta$ values are selected. First, the convergence history graphs of the fitness function over 100 iterations for different parameters are presented from Figures \ref{loss7} to \ref{loss10}. The fitness function value decreases fast over 100 iterations. 

\begin{figure}[H]
	\captionsetup{width=0.9\columnwidth}
	\centering\includegraphics[height=9cm,width=13.0cm]{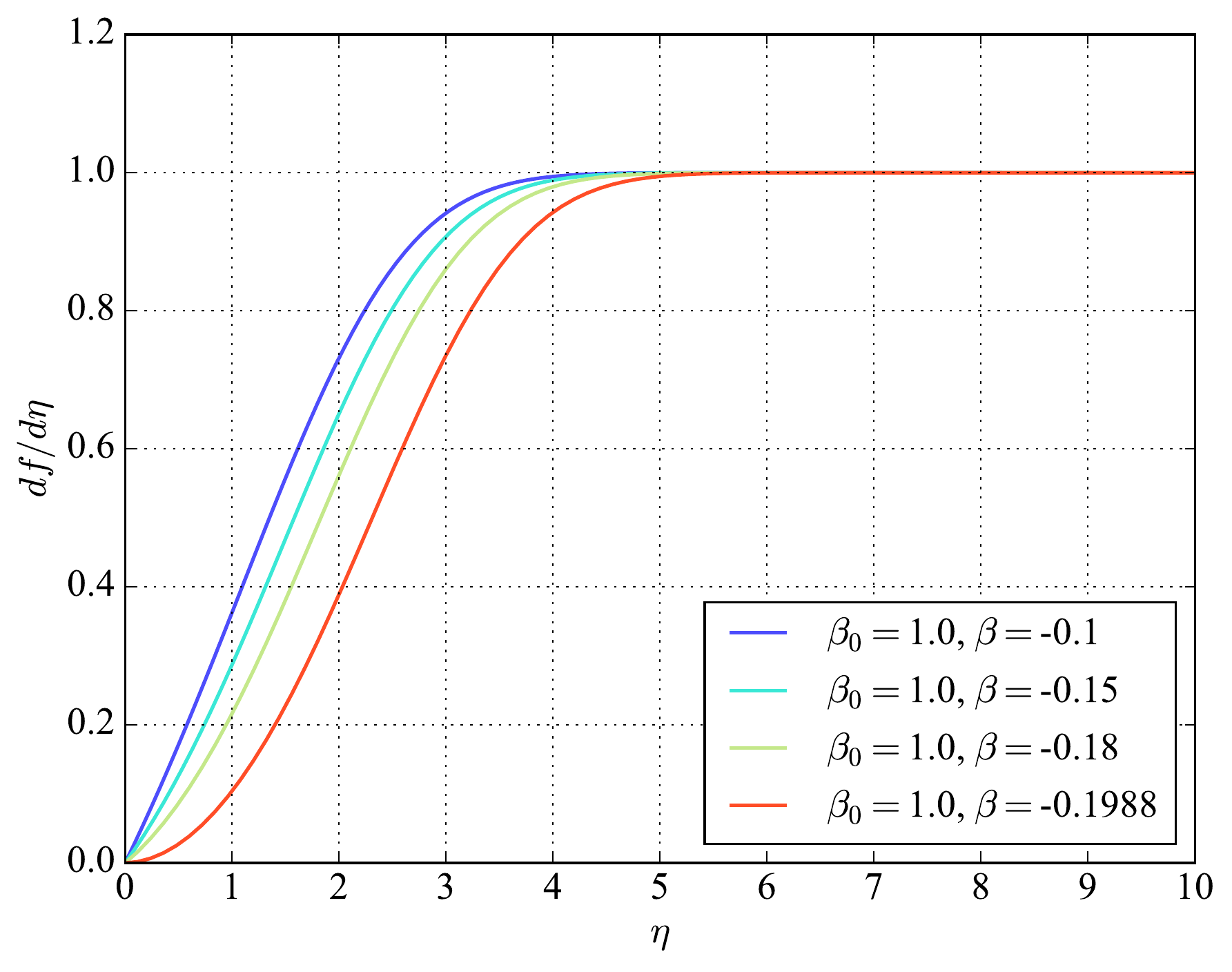}   
	\caption{ The velocity profile corresponding to different  $\beta \in \left[ -0.1988,0 \right] $ when ${ \beta  }_{ 0 }=1$}
	\label{beta3}
\end{figure}

  Fig.~\ref{beta3} illustrates the velocity profiles with the different $\beta$ values, and show thats the velocity $f'(\eta)$ increases with an increase in value of $\beta$, and the velocity profile $f'(\eta)$ approaches 1 with the $\eta$ grows. According to the results, the Runge-Kutta method combined with the Jaya algorithm has a good numerical performance to solve the Falkner-Skan equation.
For deceleration flows, we can observe that as the
parameter $\beta$ decreases, the boundary layer thickness
decreases, and eventually tends to one as the distance increases
from the initial boundary. 
\subsubsection{ Decelerating flows for $\beta_{0}$=1 and $\beta$=-0.15}
\begin{figure}[H]
	\captionsetup{width=0.9\columnwidth}
	\centering\includegraphics[height=9cm,width=13.0cm]{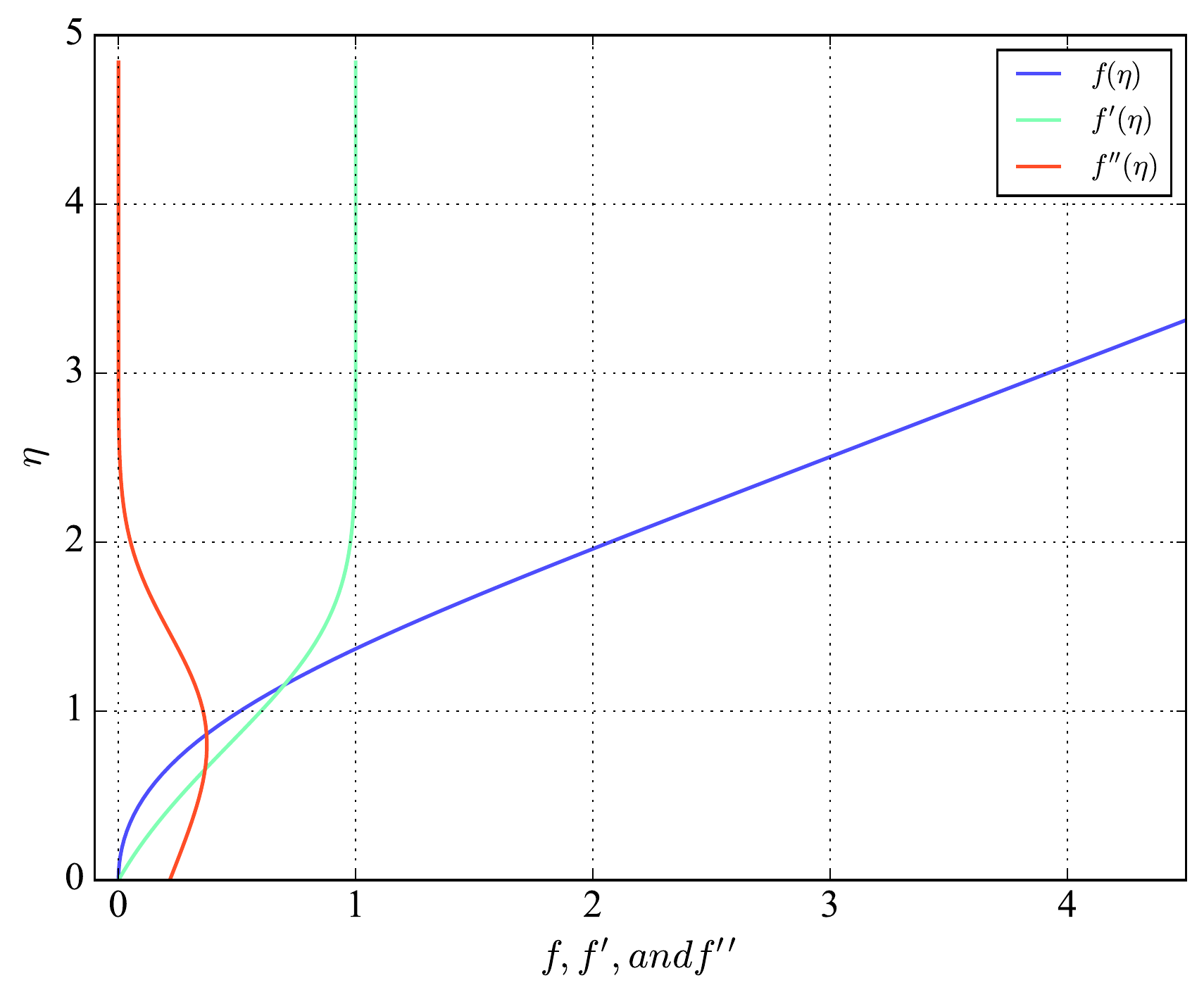}   
	\caption{ The stream function for Decelerating flow and its derivatives
		corresponding to $\beta_{0}$=1, $\beta$=-0.5}
	\label{case5}
\end{figure}

\begin{figure}[H]
	\captionsetup{width=0.9\columnwidth}
	\centering\includegraphics[height=9cm,width=13.0cm]{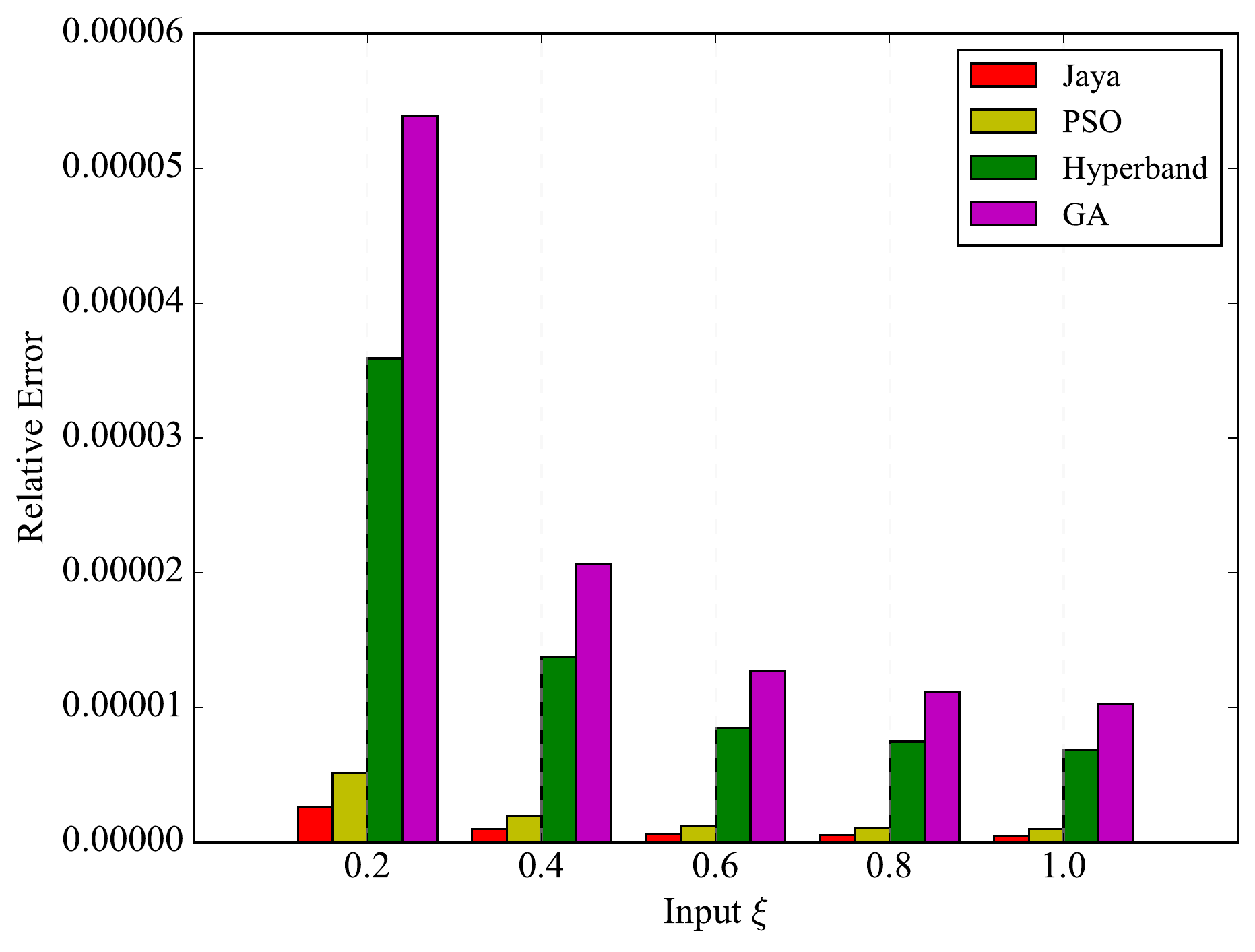}   
	\caption{Comparison of velocity of proposed results for Decelerating flow}
	\label{compare5}
\end{figure}

To be more specific, the graph of stream function $f(\eta )$, its velocity and skin friction coefficient are solved by the hybrid Jaya Runge-Kutta method illustrated in Fig.~\ref{case5}. With the increasing of $\eta$, $f^{''}$ approaches 0 and $f^{'}$ approaches 1, which agrees well with the boundary conditions in Equations \ref{bd2} and \ref{bd3}. 
Moreover, the solutions with different methods (Jaya, PSO, Hyperband, and GA) are compared with the reference solution \cite{ahmad2017stochastic} shown in Tables ~\ref{Table5}, the results obtained with Jaya method agrees pretty well with the reference solution. The absolution error for those listed optimization methods along the $\xi$ are presented Fig.~\ref{compare5}. From the graph, the Jaya algorithm obtains the stable and accurate results compared with other methods.

\begin{table}[H] 
	\captionsetup{width=0.85\columnwidth}
	\caption{Comparison of proposed results with reference solution for Decelerating flow} 
	\vspace{-0.3cm}
	\centering 
	\resizebox{0.8\columnwidth}{!}{%
		\begin{tabular}{l c c c c c} 
			\toprule 
			\toprule 
			$\xi$&$f^{'}_{ref}$&Jaya&PSO&Hyperband&GA\\
			\midrule 
    		0.1&0.25214592&0.25214504&0.25214416&0.25215822&0.25212746\\
    		\midrule
            0.2&0.57886782&0.57886633&0.57886485&0.57888861&0.57883663\\
            \midrule
            0.3&0.84962283&0.84962143&0.84962002&0.84964252&0.84959331\\
            \midrule
            0.4&0.97107618&0.97107523&0.97107427&0.97108955&0.97105613\\
            \midrule
            0.5&0.9972818&0.9972811&0.99728041&0.99729153&0.9972672\\
            \midrule
            0.6&0.99987996&0.99987935&0.99987874&0.99988845&0.99986722\\
            \midrule
            0.7&0.99999732&0.99999676&0.99999619&1.00000521&0.99998549\\
            \midrule
            0.8&0.99999974&0.99999921&0.99999868&1.0000072&0.99998856\\
            \midrule
            0.9&0.99999978&0.99999927&0.99999876&1.00000689&0.99998911\\
            \midrule
            1&0.99999979&0.9999993&0.99999881&1.00000662&0.99998954\\

			\bottomrule 
		\end{tabular}
	}
	\label{Table5} 
\end{table}

In summary, to show the ability of parameter identification with the Jaya optimizer, the $\alpha$, and $\eta_{\infty}$ results identified with different methods, which include the Runge-Kutta method combined with the Jaya algorithm, PSO algorithm, Hyperband algorithm, and GA, and the classical methods of Zhang \cite{zhang2009iterative} and Asaithambi \cite{asaithambi2005solution}, are compared in Table~\ref{Table7}. By camparing the results with different methods, we found that the results using the Runge-Kutta method combined with different algorithms have a good agreement with the results using classical methods, and the optimal results using the Runge-Kutta method combined with the Jaya algorithm are closer to the reference solutions than the other heuristic optimization methods. It can be concluded that the Runge-Kutta method combined with the Jaya algorithm is more suitable to solve Falkner-Skan equation. With the hybrid Jaya Runge-Kutta method, once the unknown parameters $\alpha$ and $\eta_{\infty}$  are determined, we can easily compute the velocity
profiles, skin friction coefficient etc. 
\begin{table}[H] 
	\captionsetup{width=1\columnwidth}
	\caption{Comparison of computed ${\alpha} $ corresponding to different $\beta$ and ${ \beta  }_{ 0 }$} 
	\vspace{-0.3cm}
	\centering 
	\resizebox{0.85\columnwidth}{!}{%
		\begin{tabular}{l c c c c c c c c c} 
			\toprule 
			\toprule 
			& & \multicolumn{3}{c}{\textbf{Jaya algorithm}}&{\textbf{Zhang}}\cite{zhang2009iterative}&{\textbf{Asaithambi}}\cite{asaithambi2005solution}&PSO&Hyperband&GA \\ 
            \cmidrule(l){3-5}
            \cmidrule(l){6-10}
			${\beta}_{0}$&${\beta}$& ${\alpha} $&$ {\eta}_{\infty}$&Residual&${\alpha} $&${\alpha} $&${\alpha} $&${\alpha} $&${\alpha} $\\ 
			\midrule 
			0.5&0&0.332057&11.856964&2.73E-24&0.33205&0.33205&0.33204&0.33142&0.33215\\ 
			\midrule
			2&1&1.311938&4.840246&1.26E-18&1.31194&1.31194&1.31185&1.31222&1.31199\\ 
			\midrule 
			1&2&1.687218& 4.547123& 2.18E-12&1.68721&1.68721&1.68772&1.68723&1.68688 \\ 
			\midrule
			1&1&1.232588& 9.078257& 5.44E-18&1.23258&1.23258&1.23254&1.23228&1.23257 \\ 
			\midrule
			1&0.5&0.927680&6.995320 &9.64E-20&0.92768&0.92768&0.92764&0.92674&0.92797 \\ 
			\midrule
			1&0&0.469600&10.746206 &6.64E-27&0.46960&0.46960&0.46957&0.47009&0.46973  \\ 
			\midrule
			1&-0.1&0.319270 &8.181430 &5.49E-21&0.31927&0.31927&0.31925&0.31815&0.31945\\ 
			\midrule
			1&-0.15&0.216361&8.975579 &3.75E-19&0.21636&0.21636&0.21636&0.21646&0.21377\\ 
			\midrule
			1&-0.18&0.128636& 11.999854&1.59E-45&0.12863&0.12863&0.12864&0.13208&0.12884\\
			\midrule
			1&-0.1988&0.005218&11.999793&1.77E-41&0.00522&0.00522&0.00559&0.00509&0.00513\\
			\bottomrule 
		\end{tabular}
	}
	\label{Table7} 
\end{table}

\section{Conclusion}
The hybrid Jaya Runge-Kutta method is presented in this paper to solve Falkner-Skan boundary value problem, which involves the identification of unknown parameters in partial differential equations. Further, this application also shows the ability of this hybrid method in solving coupled differential equations with prescribed boundary conditions. The original problems can be sensitive to the guess of initial values, with the help of Jaya algorithms, the whole scheme can yield stable and accurate results. The incompressible flow over a stretching/shrinking wedge with various wedge angles is examined. The Jaya algorithm can search the optimal parameters ${ \eta  }_{ \infty  }$ and  $\alpha$ by finding the minimal value of the fitness function according to the convergence history graph. Convergence and initial guess issues which trap those classical methods could be overcome through this simple but effective methodology.

Based on the results obtained by the hybrid Jaya Runge-Kutta method, it can be concluded that the present method can provide a reliable, effective, and accurate solution for the Falkner-Skan free boundary value problem, which includes Blasius equation, the Homann problem, the accelerating flows, and the decelerating flows. Jaya algorithm has also been verified to be effective in identification of those unknown parameters. In the future, this general of the method can be further applied into other multi-field coupled boundary layer flow problems and other time dependent partial differential equations.

 \clearpage

\bibliography{ref_FS.bib}
\end{document}